\documentclass[11pt]{article}

\usepackage[utf8]{inputenc}
\usepackage{amsmath,amsthm,amssymb,enumerate,framed}
\usepackage{color,colortab}
\usepackage[square,numbers,sort&compress]{natbib}

\usepackage{soul}

\newcommand{\Q}{\mathbb Q}
\newcommand{\N}{\mathbb{N}}
\newcommand{\Z}{\mathbb{Z}}
\newcommand{\R}{\mathbb{R}}
\newcommand{\conv}{\operatorname{conv}}

\newcommand{\rec}{\operatorname{rec}}

\newcommand{\cl}{\operatorname{cl}}

\newcommand{\cone}{\operatorname{cone}}
\newcommand{\integ}{\operatorname{integ}}

\renewcommand{\epsilon}{\varepsilon}

\newenvironment{claimproof}{\begin{trivlist} \item[] {\em Proof of Claim.}}{\hspace*{\stretch{1}} $\diamond$ \end{trivlist}}

\newtheoremstyle{mythmstyle}
	{\topsep}
	{\topsep}
	{\itshape}
	{}
	{\scshape}
	{.}
	{3pt}
	{}
\theoremstyle{mythmstyle}

\newtheorem{nn}{}[section]
\newtheorem{lemma}[nn]{Lemma}
\newtheorem{theorem}[nn]{Theorem}
\newtheorem{cor}[nn]{Corollary}
\newtheorem{prop}[nn]{Proposition}

\newtheorem{claim}[nn]{Claim}

\newtheorem{example}[nn]{Example}

\newtheorem{REMARK}[nn]{Remark}
\newenvironment{remark}{\begin{REMARK}}{\end{REMARK}}

\newcommand{\aff}{\operatorname{aff}}

\numberwithin{equation}{section}

\begin{document}

\title{The structure of the infinite models in integer programming}

\author{Amitabh Basu\footnote{Department of Applied Mathematics and Statistics, The Johns Hopkins University, {\tt basu.amitabh@jhu.edu}} \and Michele Conforti\footnote{Dipartimento di Matematica ``Tullio Levi-Civita'', Universit\`a degli Studi di Padova, Italy, {\tt conforti@math.unipd.it, disumma@math.unipd.it}}\and Marco Di Summa\footnotemark[2] \and Joseph Paat\footnote{Institute for Operations Research, Department of Mathematics, ETH Z\"urich, Switzerland, {\tt joseph.paat@ifor.math.ethz.ch}}}


\maketitle

\begin{abstract}
The infinite models in integer programming can be described as the convex hull of some points or as the intersection of halfspaces derived from valid functions. In this paper we study the relationships between these two descriptions. Our results have implications for corner polyhedra. One consequence is that nonnegative, continuous valid functions suffice to describe corner polyhedra (with or without rational data).\end{abstract}

\section{Introduction}

Let $b\in \R^n\setminus \Z^n$. The {\em mixed-integer infinite group relaxation} $M_b$ is the set of all pairs of functions $(s,y)$ with $s:\R^n\to \R_+$ and $y:\R^n\to\Z_+$ having finite support (that is, $\{r:\;s(r)>0\}$ and $\{p:\; y(p)>0\}$ are finite sets)
satisfying
$$\sum_{r\in \R^n} r s(r) + \sum_{p\in \R^n} p y(p)   \in b+\Z^n. $$

$M_b$ is a subset of the infinite-dimensional vector space $\R^{(\R^n)}\times\R^{(\R^n)}$, where $\R^{(\R^n)}$ denotes the set of functions with finite support from $\R^n$ to $\R$. Similarly, $\R^{(\R^n)}_+$ will denote the set of functions with finite support from $\R^n$ to $\R$ that are nonnegative. We will work with this vector space throughout the paper.

A tuple $(\psi, \pi, \alpha)$, where $\psi, \pi: \R^n \to \R$ and $\alpha \in \R$, is a {\em valid tuple for $M_b$} if
\begin{equation}\label{eq:valid-Mb}\sum_{r\in \R^n} \psi(r)s(r) + \sum_{p\in \R^n} \pi(p)y(p)\ge \alpha \mbox{ for every }(s,y) \in M_b. \end{equation}
A valid  tuple  $(\psi,\pi,\alpha)$ for $M_b$ is {\em minimal} if there does not exist a pair of functions $(\psi', \pi')$ different from $(\psi, \pi)$, with $(\psi', \pi')\le(\psi, \pi)$, such that $(\psi', \pi', \alpha)$ is a valid tuple for $M_b$. Since for $\lambda>0$ the inequalities \eqref{eq:valid-Mb} associated with $(\psi, \pi, \alpha)$ and $(\lambda\psi, \lambda\pi, \lambda\alpha)$ are  equivalent, from now on we assume $\alpha \in \{-1,0,1\}$. \medskip

The set of functions $y:\R^n\to \Z_+$ such that $(0,y) \in M_b$ will be called the {\em pure integer infinite group relaxation} $I_b$. In other words, $I_b =\{ y \in \R^{(\R^n)}: (0,y) \in M_b\}$. When convenient we will see $I_b$ as a subset of $M_b$. A tuple $(\pi, \alpha)$, where $\pi: \R^n \to \R$ and $\alpha \in \R$, is called a {\em valid tuple for $I_b$} if
\begin{equation}\label{eq:valid-Ib}
\sum_{p\in \R^n} \pi(p)y(p)\ge \alpha\mbox{ for every }y \in I_b.
\end{equation}
Again, we will assume $\alpha\in\{-1,0,1\}$. We say that a valid tuple $(\pi,\alpha)$ for $I_b$ is {\em minimal} if there does not exist a function $\pi'$ different from $\pi$, with $\pi'\le\pi$, such that $(\pi', \alpha)$ is a valid tuple for $I_b$.
\smallskip

Models $M_b$ and $I_b$ were defined by Gomory and Johnson in a series of papers~\cite{gom,infinite,infinite2,johnson1974group} as a template to generate valid inequalities for integer programs. They have been the focus of extensive research, as summarized, e.g., in \cite{basu2015geometric,basu2016light,basu2016light2}, \cite[Chapter 6]{conforti2014integer}. The main idea is the following. Given a mixed-integer program in the form $Cz = d, z\geq 0$ with integrality constraints on a subset of the $z$ variables, one can obtain the following relaxation from an optimal simplex tableau:

$$x + b = \sum_{r\in R} r s(r) + \sum_{p\in P} p y(p),$$
where $x$ is a subset of $n$ integer constrained basic variables, $R, P \subseteq \R^n$ are the sets of nonbasic columns corresponding to the continuous and integer constrained nonbasic variables $s(r)$ and $y(p)$ respectively, and $(x,s,y) = (-b,0,0)$ is the current linear-programming (LP) solution. If one relaxes the nonnegativity constraint on $x$, the relaxation can be seen as a finite dimensional face of $\conv(M_b)$ (the convex hull of $M_b$), obtained by setting $s(r)$ and $y(p)$ to $0$ for $r \not\in R$ and $p\not\in P$. Similarly, for a pure integer program, one obtains a finite dimensional face of $\conv(I_b)$. The main point is that any valid tuple (for $M_b$ or $I_b$), when restricted to the appropriate finite dimensional space, gives  a valid inequality which can be used as a cutting plane for the initial LP relaxation. 

In the pure integer case, when all the data in the problem is rational, i.e., $P\cup\{b\} \subseteq \Q^n$, Gomory termed such finite dimensional faces of $\conv(I_b)$ as {\em corner polyhedra}. One of our insights in this paper is that even for non-rational data, the finite dimensional faces of $\conv(I_b)$ are rational polyhedra (see detailed discussion below). In anticipation, we will refer to all finite dimensional faces of $\conv(I_b)$ (with rational or non-rational $P, b$) as {\em corner polyhedra}. 


\subsection{Results and their implications.} We first study the infinite dimensional objects $M_b$ and $I_b$, and then derive some consequences for the finite dimensional faces of $\conv(M_b)$ and $\conv(I_b)$. For later reference, we define a {\em canonical face of  $\conv (M_b)$} as  $\conv (M_b)\cap \{(s,y):  s(r)=0\; \forall r\in \R^n\setminus R,\, y(r)=0\; \forall r\in \R^n\setminus P\}$ for some $R,P\subseteq \R^n$.  Similarly, a {\em canonical face of  $\conv (I_b)$} is  $\conv (I_b)\cap \{y:  y(r)=0\;\forall r\in \R^n\setminus P\}$ for some $P\subseteq \R^n$. These will be called {\em finite canonical faces} if $R, P$ are finite subsets of $\R^n$, which will be the focus of our investigations. 

\paragraph{Structure of the infinite dimensional models.} Regarding the infinite dimensional objects $M_b$ and $I_b$, one would expect that the intersection of $\R^{(\R^n)}_+\times\R^{(\R^n)}_+$ and all halfspaces in $\R^{(\R^n)}\times\R^{(\R^n)}$ defined by valid tuples for $M_b$ would be equal to $\conv(M_b)$, where $\conv(\cdot)$ denotes the convex hull operator. However, this is not true: one of our main results (Theorem \ref{thm:closure-char}) shows that this intersection is the closure of $\conv(M_b)$ under a norm topology on  $\R^{(\R^n)}\times\R^{(\R^n)}$ that was first defined by Basu et al. \cite{bccz}. We then show that the closure of $\conv(M_b)$ coincides with $\conv (M_b)+(\R_+^{(\R^n)}\times\R_+^{(\R^n)})$, which is a strict superset of $\conv(M_b)$ (Remark~\ref{rem:Q_b-strict}). The same set of results holds for $\conv(I_b)$ (Theorem~\ref{thm:closure-char-Ib}, Remark~\ref{rem:Gb-strict} and Example~\ref{ex:pure-integer}). We also obtain precise characterizations of the affine hulls of $M_b$ and $I_b$ (Propositions~\ref{prop:aff-hull-M_b} and~\ref{prop:affine-hull} respectively); our characterization of the equations defining the affine hull of $\conv(I_b)$ extends a result in~\cite{basu2016light}.

Crucial to the above results is the characterization of what we call {\em liftable tuples}. We say that a valid tuple $(\pi,\alpha)$ for $I_b$ is {\em liftable} if there exists  $\psi:\R^n\to \R$ such that $(\psi,\pi,\alpha)$ is a valid tuple for $M_b$. Minimal valid tuples $(\pi,\alpha)$ that are liftable are a strict subset of minimal valid tuples, as we show that such $\pi$ have to be nonnegative and Lipschitz continuous (Proposition \ref{prop:liftable} and Remark~\ref{rem:sup-lipschitz}).

\paragraph{Restricting to well-behaved valid tuples.} Most of the prior literature on valid tuples $(\pi,\alpha)$ for $I_b$ proceeds under the assumption that $\pi$ is nonnegative (in fact, Gomory and Johnson included the assumption $\pi\geq 0$ in their original definition of valid tuple for $I_b$). This assumption is restrictive as there are valid tuples not satisfying $\pi\ge0$~\cite{basu2016light}, and the assumption has been investigated in more recent work on $I_b$~\cite{basu2016light,basu2016light2} (work on {\em generalizations} of $I_b$ also allows for negative values in the valid tuples~\cite{yildiz2015integer,kilincc2015sufficient,kilincc2016sublinear,kilincc2015minimal}). The standard justification behind the $\pi\ge0$ assumption is the fact that valid tuples are nonnegative on the rational vectors; see, e.g., the discussion in Section 2.1.2 in~\cite{basu2016light}. Since in practice we are interested in finite canonical faces of $\conv(I_b)$ defined by $P\cup\{b\}\subset \Q^n$, such an assumption seems reasonable. However, no mathematical evidence exists in the literature that a complete inequality description of the finite canonical faces (even with rational data) can be obtained from the nonnegative valid tuples only.\footnote{Such results are obtainable in the one-dimensional ($n=1$) rational case by elementary means such as interpolation (using lifting techniques for finite dimensional rational polyhedra and results like~\cite[Theorem 8.3]{basu2016light2}; see also the recent work in~\cite{koppe2018all}). We are unaware of a way to establish these results for general $n\geq 2$ without using the technology developed in this paper.}

In this paper we give concrete mathematical evidence that {\em one can restrict to nonnegative valid tuples without any loss of generality}. Being able to restrict to nonnegative valid tuples has the added advantage that nonnegative minimal valid tuples form a {\em compact, convex set} under the natural product topology on functions. Thus, one approach to understanding valid tuples is to understand the extreme points of this compact convex set, which are termed {\em extreme functions/tuples} in the literature. While this nonnegativity assumption was standard for the area, our results discussed below about nonnegative valid tuples now give a rigorous justification for this.
\medskip

First, we show that for every minimal valid tuple $(\pi,\alpha)$, there exist a {\em unique} $\theta:\R^n\to\R$ and $ a \in \R$ such that  both $(\theta,  a)$, $(-\theta, - a)$ are minimal valid tuples and the minimal valid tuple $(\pi',\alpha')=(\pi-\theta,\alpha-a)$ satisfies $\pi'\ge0$ (Theorem \ref{th:nonneg}). In other words, every minimal valid tuple is equivalent to a (unique) nonnegative minimal valid tuple, modulo an equation for the affine hull of $I_b$. Since every valid tuple is dominated by a minimal valid tuple, this shows that every valid tuple is dominated by a nonnegative minimal valid tuple, modulo the affine hull of $I_b$. This settles an open question in \cite[Open Question 2.5]{basu2016light}.


Second, we show that any finite canonical face $F$ of $\conv(I_b)$ is a rational polyhedron, even if $P\cup\{b\}$ contains non-rational vectors (Theorem \ref{thm:rec-cone-finite}). This justifies our use of the term {\em corner polyhedron} to refer to {\em any} finite canonical face, extending Gomory's original use of {\em corner polyhedron} which applied only to finite canonical faces with rational data. Moreover, we prove that $F$ is given by the intersection of $\aff(F)$ and the inequalities obtained by the restriction of minimal, liftable tuples for $I_b$ (Theorem~\ref{thm:conv-Ib-nonneg}). Since liftable tuples are always nonnegative and Lipschitz continuous, this reinforces the claim that one can restrict attention to only nonnegative minimal tuples. Theorem~\ref{thm:conv-Ib-nonneg} also has the infinite dimensional interpretation that $\conv(I_b)$ is given by the intersection of $\R^{(\R^n)}_+$ with the affine hull $\aff(I_b)$ and the halfspaces given by liftable tuples (Corollary~\ref{cor:conv-Ib-nonneg}). Combined with the structural results for $I_b$ mentioned above, this also implies that $\conv(I_b) = \cl(\conv(I_b)) \cap \aff(I_b) = (\conv(I_b)+\R^{(\R^n)}_+) \cap \aff(I_b)$, where $\cl(\cdot)$ denotes the closure operator with respect to the norm topology discussed above.

Finally, we strengthen the above result for corner polyhedra with rational data, i.e., $P \subseteq \Q^n$. Theorem~\ref{thm:conv-Ib-nonneg} only gives the guarantee that when we intersect {\em all} liftable tuples (restricted to the appropriate finite dimensional space), we obtain the corner polyhedron. Theorem~\ref{thm:conv-Ib-nonneg} does not rule out the possibility that for a given corner polyhedron and a particular facet defining inequality for it, this inequality can only be obtained as a ``limit" of minimal liftable tuples; i.e., there is no single liftable tuple that dominates the given facet defining inequality. We show that if we have {\em rational} data, i.e., $P \subseteq \Q^n$, then this issue does not occur (Theorem~\ref{thm:facet-restriction}). In particular, we show that if a corner polyhedron is nonempty, then any nontrivial valid inequality for it (i.e., an inequality that is not implied by the nonnegativity constraints on the $y$ variables) is dominated by the restriction of an inequality given by some minimal liftable tuple. Consequently, for any finite $P \subseteq \Q^n$, the associated corner polyhedron can be given a complete, finite description by the restriction of some minimal, liftable tuples and the nonnegativity constraints on the $y$ variables.


\medskip

Literature on valid tuples contains constructions of families of extreme valid tuples $(\pi,\alpha)$ such that $\pi$ is discontinuous~\cite{Miller-Li-Richard2008,koppe2015electronic,letchford2002strengthening,dey1,twoStepMIR,hildebrand-thesis} (or continuous but not Lipschitz continuous~\cite{koppe2015electronic}). Our result above shows that such functions may be disregarded, if one is interested in valid inequalities or facets of rational corner polyhedra. Our results show that such extreme tuples are redundant within the set of valid tuples, as far as rational corner polyhedra are concerned. Moreover, even for corner polyhedra with non-rational data, one can restrict attention to liftable tuples (and therefore nonnegative and Lipschitz continuous valid tuples) as long as one is interested in a (not necessarily finite) halfspace description of the corner polyhedron. We also note that related results on sufficiency of valid tuples were obtained in~\cite{zambelli,conforti2010equivalence,cornuejols2015sufficiency,kilincc2015sufficient,kilincc2016sublinear}.

\paragraph{Geometry of corner polyhedra.} The above results are derived out of a detailed study of the geometry of corner polyhedra, which is interesting in its own right. As mentioned above, Theorem \ref{thm:rec-cone-finite} shows that the finite canonical faces of $\conv(I_b)$ are always rational polyhedra, justifying the use of ``corner polyhedra" even in the presence of non-rational data. Theorem \ref{thm:rec-cone-finite} also characterizes the recession cone of such a face; it is simply the intersection of the nonnegative orthant and the linear space parallel to the affine hull of $\conv(I_b)$, restricted to the finite dimensional space containing the face.
This proves to be a crucial insight. Theorem~\ref{thm:rational-P} sharpens these results to give a tight characterization of corner polyhedra with rational data. In particular, it shows that a corner polyhedron is defined using rational data {\em if and only if} its recession cone is the nonnegative orthant, which happens {\em if and only if} the corner polyhedron is full-dimensional. Both of these theorems rely on the structure of the infinite models unveiled in this paper. It is interesting to note that several important properties of finite dimensional corner polyhedra -- which are most relevant for integer programming -- are revealed through a detailed study of the infinite dimensional models.

\paragraph{Topological pathologies in finite canonical faces of $M_b$.} While the finite canonical faces of $\conv(I_b)$ are always rational polyhedra, the finite canonical faces of $\conv (M_b)$ are not as well-behaved: there exist finite canonical faces of $M_b$ that are not closed (in the standard finite dimensional topology) -- see Example \ref{ex:not-closed}. 

\medskip

The remainder of the paper is dedicated to establishing the results discussed above.

\section{The structure of $\conv(M_b)$ and $\conv(I_b)$.}\label{sec:structure}

We start with a well-known fact about minimal valid tuples.

\begin{remark}\label{rem:Zorn}
An application of Zorn's lemma (see, e.g., \cite[Proposition A.1]{basu-paat-lifting}) shows that, given a valid tuple $(\psi,\pi,\alpha)$ for $M_b$, there exists a minimal valid tuple $(\psi',\pi',\alpha)$ for $M_b$ with $\psi'\le\psi$ and $\pi'\le\pi$. Similarly, given a valid tuple $(\pi,\alpha)$ for $I_b$, there exists a minimal valid tuple $(\pi',\alpha)$ for $I_b$ with $\pi'\le\pi$. We will use this throughout the paper.
\end{remark}

Given a tuple $(\psi, \pi,\alpha)$, we define
$$H_{\psi,\pi,\alpha} := \bigg\{(s,y)\in\R^{(\R^n)} \times \R^{(\R^n)}: \sum_{r\in \R^n} \psi(r)s(r) + \sum_{p\in \R^n} \pi(p)y(p)\ge \alpha\bigg\}.$$

A valid tuple $(\psi,\pi,\alpha)$ for $M_b$ is {\em trivial} if $\R^{(\R^n)}_+ \times \R^{(\R^n)}_+\subseteq H_{\psi,\pi,\alpha}$. This happens if and only if $\psi\ge 0$, $\pi \ge 0$ and $\alpha\in\{0,-1\}$. Similarly, a valid tuple $(\pi,\alpha)$ for $I_b$ is {\em trivial} if $\pi\ge 0$ and $\alpha\in\{0,-1\}$. 
\smallskip

A function $\phi:\R^n\to\R$ is {\em subadditive} if $\phi(r_1)+\phi(r_2)\ge\phi(r_1+r_2)$ for every $r_1,r_2\in\R^n$,
and is {\em positively homogenous} if $\phi(\lambda r)=\lambda\phi(r)$ for every $r\in\R^n$ and $\lambda\ge0$.
If $\phi$ is subadditive and positively homogenous, then $\phi$ is called {\em sublinear}. The following proposition is well-known and its proof can be found in the Appendix.

\begin{prop}\label{prop:psi_sublinear}
Let $(\psi, \pi, \alpha)$ be a minimal valid tuple for $M_b$. Then $\psi$ is sublinear and $\pi\leq \psi$.
\end{prop}

\begin{lemma}\label{lem:psi-using-sup}
Suppose $\pi:\R^n\to\R$ is subadditive and $\sup_{\epsilon > 0}\frac{\pi(\epsilon r)}{\epsilon} < \infty$ for all $r \in \R^n$. Define $\psi:\R^n\to \R$ by $\psi(r) := \sup_{\epsilon > 0}\frac{\pi(\epsilon r)}{\epsilon}$. Then $\psi$ is sublinear and $\pi\le\psi$.
\end{lemma}

\proof{Proof.}
Since $\pi$ is subadditive, $\psi$ is readily checked to be subadditive as well. The fact that $\pi\le\psi$ follows by taking $\epsilon=1$. Finally, positive homogeneity of $\psi$ follows from the definition of $\psi$.
\endproof


\begin{theorem}\label{thm:non-dominated} Let $\psi:\R^n\to\R$, $\pi:\R^n\to\R$ be any functions, and $\alpha\in\{-1,0,1\}$. Then $(\psi, \pi, \alpha)$ is a nontrivial minimal valid tuple for $M_b$ if and only if all of the following hold:
\begin{enumerate}[(a)]
\item $\pi$ is subadditive;
\item $\psi(r) = \sup_{\epsilon > 0}\frac{\pi(\epsilon r)}{\epsilon} = \lim_{\epsilon \to 0^+}\frac{\pi(\epsilon r)}{\epsilon} = \lim\sup_{\epsilon \to 0^+}\frac{\pi(\epsilon r)}{\epsilon}$ for every $r\in \R^n$;
    \item $\pi$ is Lipschitz continuous with Lipschitz constant $L:=\max_{\|r\|=1}\psi(r)$;
\item $\pi \geq 0$, $\pi(z)=0$ for every $z\in \Z^n$, and $\alpha = 1$;
\item (symmetry condition) $\pi$ satisfies $\pi(r) + \pi(b-r) = 1$ for all $r\in \R^n$.
\end{enumerate}
\end{theorem}

The above theorem, except for the Lipschitz continuity and nonnegativity of $\pi$, follows from a result of Y\i ld\i z and Cornu\'ejols \cite[Theorem 37]{yildiz2015integer} together with the characterization of the nontrivial minimal valid tuples for $I_b$ due to Gomory and Johnson (see, e.g., \cite[Theorem 6.22]{conforti2014integer}). See also the result of Johnson for valid tuples for $M_b$~\cite[Theorem 6.34]{conforti2014integer}. We provide a self-contained proof of Theorem~\ref{thm:non-dominated} in the Appendix.

\begin{cor}\label{cor:dominated_minimal}
Let $(\pi,\alpha)$ be a nontrivial minimal valid tuple for $I_b$ such that $\sup_{\epsilon>0}\frac{\pi(\epsilon r)}{\epsilon}<\infty$ for every $r\in\R^n$. Define $\psi:\R^n\to \R$ by $\psi(r):=\sup_{\epsilon>0}\frac{\pi(\epsilon r)}{\epsilon}$. Then $(\psi,\pi,\alpha)$ is a nontrivial minimal valid tuple for $M_b$.

Conversely, if $(\psi,\pi,\alpha)$ is a nontrivial minimal valid tuple for $M_b$, then $(\pi,\alpha)$ is a nontrivial minimal valid tuple for $I_b$.
\end{cor}

\proof{Proof.}
Since $(\pi,\alpha)$ is minimal, the same argument as in the proof of Proposition \ref{prop:psi_sublinear} shows that $\pi$ is subadditive. Let $\psi$ be defined as above.
Following the proof of Theorem \ref{thm:non-dominated} it can be checked that minimality and nontriviality of $(\pi,\alpha)$ suffice  to show that $(\psi,\pi,\alpha)$ satisfies $(a)$--$(e)$ of Theorem \ref{thm:non-dominated}, and therefore $(\psi,\pi,\alpha)$ is a nontrivial minimal valid tuple for $M_b$.

For the converse, we use a theorem of Gomory and Johnson (see, e.g., \cite[Theorem 6.22]{conforti2014integer}) stating that
if $(\pi,1)$ is a nontrivial valid tuple with $\pi\ge0$, then $(\pi,1)$ is minimal if and only if $\pi$ is subadditive, $\pi(z)=0$ for every $z\in\Z^n$, and $\pi$ satisfies the symmetry condition.
Let $(\psi,\pi,\alpha)$ be a nontrivial minimal valid tuple for $M_b$.
By Theorem \ref{thm:non-dominated}, $\pi\ge0$, $\alpha=1$, $\pi$ is subadditive, $\pi(z)=0$ for every $z\in\Z^n$, and $\pi$ satisfies the symmetry condition. Therefore, by the above theorem, $(\pi,\alpha)$ is a nontrivial minimal valid tuple for $I_b$.
\endproof

A valid  tuple $(\pi, \alpha)$ for $I_b$ is called {\em liftable} if there exists a function $\psi:\R^n\to \R$ such that $(\psi, \pi, \alpha)$ is a valid tuple for $M_b$.



\begin{prop}\label{prop:liftable}
Let $(\pi,\alpha)$ be a nontrivial valid tuple for $I_b$. Then $(\pi,\alpha)$ is liftable if and only if there exists a minimal valid tuple $(\pi',\alpha)$ such that $\pi'\le\pi$ and $\sup_{\epsilon > 0}\frac{\pi'(\epsilon r)}{\epsilon}<\infty$ for every $r \in \R^n$. In this case, defining $\psi(r)=\sup_{\epsilon > 0}\frac{\pi'(\epsilon r)}{\epsilon}$, we have that $(\psi,\pi',\alpha)$ is a nontrivial valid tuple  for $M_b$ which is minimal .
\end{prop}

\proof{Proof.}
If $(\pi,\alpha)$ is nontrivial and liftable, then there exists $\psi$ such that $(\psi, \pi, \alpha)$ is a valid tuple for $M_b$. Let $(\psi', \pi', \alpha)$ be a minimal valid tuple with $\psi'\le\psi$ and $\pi'\le\pi$. Since $(\pi,\alpha)$ is nontrivial, so is $(\psi', \pi', \alpha)$. By Theorem~\ref{thm:non-dominated}, $\sup_{\epsilon > 0}\frac{\pi'(\epsilon r)}{\epsilon} < \infty$ for every $r \in \R^n$. By Corollary \ref{cor:dominated_minimal}, $(\pi',\alpha)$ is minimal.

Conversely, let $(\pi,\alpha)$ be a nontrivial valid tuple for $I_b$, and let $\pi'\le\pi$ be such that $(\pi',\alpha)$ is minimal (and nontrivial) and  $\psi(r):=\sup_{\epsilon > 0}\frac{\pi'(\epsilon r)}{\epsilon}$ is finite for every $r \in \R^n$. By Corollary \ref{cor:dominated_minimal}, $(\psi,\pi',\alpha)$ is a nontrivial minimal valid tuple for $M_b$, and therefore $(\pi',\alpha)$ is liftable. Since $\pi\ge\pi'$, $(\psi, \pi, \alpha)$ is a valid tuple for $M_b$ and therefore $(\pi,\alpha)$ is liftable as well.
\endproof




\begin{remark}\label{rem:sup-lipschitz}
Let $(\pi,\alpha)$ be a nontrivial valid tuple for $I_b$ that is minimal and liftable. It follows from Proposition \ref{prop:liftable} (with $\pi'=\pi$) that $\psi(r):=\sup_{\epsilon > 0}\frac{\pi(\epsilon r)}{\epsilon}$ is finite for all $r\in\R^n$ and $(\psi,\pi,\alpha)$  nontrivial valid tuple for $M_b$ that is minimal. Therefore by Theorem~\ref{thm:non-dominated}, $\pi$ is Lipschitz continuous and $\pi\ge0$.

There are nontrivial minimal valid tuples $(\pi,\alpha)$ for $I_b$ for which $\pi$ is not continuous, or $\pi$ is continuous but not Lipschitz continuous, see the construction in~\cite[Section 5]{koppe2015electronic}. There are also nontrivial minimal valid tuples $(\pi,\alpha)$ for $I_b$ with $\pi\not\ge0$. By the above results, none of these minimal tuples are liftable.
\end{remark}

\subsection{The closure of $\conv(M_b)$.}

\begin{lemma}\label{lem:equal-sets}
The following sets coincide:
\begin{enumerate}[(a)]
\item $\big(\R^{(\R^n)}_+ \times \R^{(\R^n)}_+\big) \cap \bigcap\{H_{\psi,\pi,\alpha} : (\psi,\pi,\alpha) \mbox{ valid tuple} \}$
\item $\big(\R^{(\R^n)}_+ \times \R^{(\R^n)}_+\big) \cap \bigcap\{H_{\psi,\pi,\alpha} : (\psi,\pi,\alpha) \mbox{ nontrivial valid tuple} \}$
\item $\big(\R^{(\R^n)}_+ \times \R^{(\R^n)}_+\big) \cap \bigcap\{H_{\psi,\pi,\alpha} : (\psi,\pi,\alpha) \mbox{ nontrivial minimal valid tuple} \}$
\item $\big(\R^{(\R^n)}_+ \times \R^{(\R^n)}_+\big) \cap \bigcap\{H_{\psi,\pi,\alpha} : (\psi,\pi,\alpha) \mbox{ nontrivial minimal valid tuple},\, \psi,\pi\ge0,\,\alpha=1 \}$
\end{enumerate}
\end{lemma}

\proof{Proof.}
The equivalence of $(a)$ and $(b)$ follows from the definition of nontrivial valid tuple.
The sets $(b)$ and $(c)$ coincide by Remark \ref{rem:Zorn}. Finally, Theorem \ref{thm:non-dominated} shows that $(c)$ is equal to $(d)$.
\endproof

From now on, we denote by $Q_b$ the set(s) of Lemma \ref{lem:equal-sets}.\medskip

While $\conv(M_b)\subseteq Q_b$, this containment is strict, as shown in Remark \ref{rem:Q_b-strict}.
However, Theorem \ref{thm:closure-char} below proves that, under an appropriate topology, the closure of $\conv(M_b)$ is exactly $Q_b$.
In order to show this result, we need the following lemma, that may be of independent interest.

\begin{lemma}\label{lemma:closed}
If $C\subseteq\R^n_+$ is closed, then so is $\conv(C)+\R^n_+$.
\end{lemma}

\proof{Proof.}
Let $(x_i)_{i\in \N}$ be a sequence of points in $\conv(C)+\R^n_+$ that converges to some $\bar x\in\R^n$. We need to show that $\bar x\in\conv(C)+\R^n_+$.

By Carath\'eodory's theorem, for every $i\in\N$ we can write
\begin{equation}\label{eq:comb}
x_i=\sum_{t=1}^{n+1}\lambda_i^tx_i^t+r_i,
\end{equation}
where $x_i^t\in C$ for all $t$, $\lambda_i^t\ge0$ for all $t$, $\sum_t \lambda_i^t=1$, and $r_i\in\R^n_+$.

Since $C$ is a closed set and the interval $[0,1]$ is compact, by repeatedly taking subsequences of the original sequence $(x_i)_{i\in\N}$, we assume that for every $t=1,\dots,n+1$ the following two conditions hold:
\begin{enumerate}[(a)]
\item either the sequence $(x_i^t)_{i\in\N}$ is unbounded or it converges to some $\bar x^t\in C$;
\item the sequence $(\lambda_i^t)_{i\in\N}$ converges to some number $\bar\lambda^t\in[0,1]$.
\end{enumerate}
Note that $\sum_{t=1}^{n+1}\bar\lambda^t=1$.

Let $T_1\subseteq\{1,\dots,n+1\}$ be the set of indices such that the sequence $(x_i^t)_{i\in\N}$ converges to $\bar x^t$, and let $T_2=\{1,\dots,n+1\}\setminus T_1$. For $i\in\N$ we rewrite \eqref{eq:comb} as
\begin{equation}\label{eq:comb2}
x_i-\sum_{t\in T_1}\lambda_i^tx_i^t=\sum_{t\in T_2}\lambda_i^tx_i^t+r_i.
\end{equation}
Since the left-hand side of \eqref{eq:comb2} converges to
\begin{equation}\label{eq:comb3}
\bar r:=\bar x-\sum_{t\in T_1}\bar\lambda^t\bar x^t,
\end{equation}
the right-hand side must also converge to $\bar r$.
Note that $\bar r\in\R^n_+$, as the right-hand side of \eqref{eq:comb2} is a nonnegative vector for all $i\in\N$.
Furthermore, $\bar\lambda^t=0$ for every $t\in T_2$, otherwise the right-hand side of  \eqref{eq:comb2} would not converge.
This implies that $\sum_{t\in T_1}\bar\lambda^t=1$ and thus equation \eqref{eq:comb3} proves that $\bar x\in\conv(C)+\R^n_+$.
\endproof

Define the following norm on $\R^{(\R^n)} \times \R^{(\R^n)}$, which was first introduced in \cite{bccz}:
$$|(s,y)|_* := |s(0)| + \sum_{r\in \R^n} \|r\||s(r)| + |y(0)| + \sum_{p\in \R^n} \|p\||y(p)|.$$

For any two functions $\psi:\R^n \to \R$, $\pi : \R^n \to \R$, we define a linear functional $F_{\psi, \pi}$ on the space $\R^{(\R^n)}\times \R^{(\R^n)}$ as follows:
$$F_{\psi,\pi}(s,y) := \sum_{r\in \R^n} \psi(r)s(r) + \sum_{p\in \R^n} \pi(p)y(p).$$

\begin{lemma}\label{lem:continuous}
Under the $|(\cdot, \cdot)|_*$ norm, the linear functional $F_{\psi,\pi}$ is continuous if $(\psi,\pi,1)$ is a nontrivial minimal valid tuple for $M_b$.
\end{lemma}

\proof{Proof.}
Since $(\psi,\pi,1)$ is a nontrivial minimal valid tuple for $M_b$, conditions $(a)$--$(e)$ of Theorem \ref{thm:non-dominated} are satisfied.
In order to show that $F_{\psi,\pi}$ is continuous, it is equivalent to show that $F_{\psi,\pi}$ is bounded, i.e., there exists a number $M$ such that $|F_{\psi,\pi}(s,y)| \leq M$ for all $(s,y)$ satisfying $|(s,y)|_* = 1$ (see Conway~\cite[Chapter III, Proposition 2.1]{conway2013course}).

We claim that $M$ can be chosen to be $\max_{\|r\|=1}\psi(r)$. (This maximum exists because, by condition $(b)$ in Theorem~\ref{thm:non-dominated}, $\psi$ is sublinear and therefore continuous on $\R^n$.) Consider $(s,y)$ such that $|(s,y)|_* = 1$. Using $\pi \leq \psi$ (Proposition~\ref{prop:psi_sublinear}), $\pi, \psi\ge0$ and $\psi(0) = \pi(0)=0$ (Theorem \ref{thm:non-dominated}), we have
\[\begin{split}
|F_{\psi, \pi}(s,y)| & =\textstyle \left|\sum_{r\in \R^n\setminus \{0\}} \psi(r)s(r) + \sum_{p\in \R^n\setminus \{0\}} \pi(p)y(p)\right|\\
& \textstyle\le \sum_{r\in \R^n\setminus \{0\}} \psi(r)|s(r)| + \sum_{p\in \R^n\setminus \{0\}} \psi(p)|y(p)| \\
& =\textstyle\sum_{r\in \R^n\setminus \{0\}} \psi\left(\frac{r}{\|r\|}\right)\|r\||s(r)| + \sum_{p\in \R^n\setminus\{0\}} \psi\left(\frac{p}{\|p\|}\right)\|p\||y(p)|\\
& \leq \textstyle M \left(\sum_{r\in \R^n\setminus \{0\}} \|r\||s(r)| + \sum_{p\in \R^n\setminus \{0\}} \|p\||y(p)|\right) \\
& \leq \textstyle M \left(|s(0)| + \sum_{r\in \R^n} \|r\||s(r)| +|y(0)|+ \sum_{p\in \R^n} \|p\||y(p)|\right) \\
& = M|(s,y)|_* = M. 
\end{split}\]
\endproof

\begin{lemma}\label{lem:finite-support-closed} Under the $|(\cdot, \cdot)|_*$ norm, the linear functional $F_{\psi,\pi}$ is continuous if $\psi$ and $\pi$ have finite supports.
\end{lemma}

\proof{Proof.} Let $R, P\subseteq\R^n$ be the supports of $\psi, \pi$ respectively. We assume $R\cup P\ne\emptyset$, otherwise the continuity of $F_{\psi,\pi}$ is obvious.
Define
$$N := \max\left\{1,\max_{r \in R\setminus\{0\}}\frac{1}{\|r\|}, \, \max_{p\in P\setminus\{0\}}\frac{1}{\|p\|}\right\},\quad
L := \max\left\{\max_{r\in R}|\psi(r)|,\,\max_{p\in P}|\pi(p)|\right\},$$
and $M := N\cdot L$. One now checks that
\[\begin{split}
|F_{\psi,\pi}(s,y)|& =\textstyle \left|\sum_{r\in R}\psi(r)s(r) + \sum_{p\in P}\pi(p)y(p)\right|
\\ &\leq \textstyle L\left(\sum_{r\in R}|s(r)| + \sum_{p\in P}|y(p)|\right) \\
& \leq \textstyle L N \left(|s(0)| + \sum_{r\in R\setminus\{0\}} \|r\||s(r)| + |y(0)| + \sum_{p\in P\setminus\{0\}} \|p\||y(p)|\right)\\
& = M |(s,y)|_*.
\end{split}\]
This shows that $F_{\psi,\pi}$ is a bounded linear functional, and hence continuous.
\endproof

\begin{lemma}\label{lem:Qb-closed}
Under the topology induced by $|(\cdot, \cdot)|_*$, the set $Q_b$ is closed.
\end{lemma}

\proof{Proof.}
Since $\R^{(\R^n)}_+ \times \R^{(\R^n)}_+$ is the intersection of a family of halfspaces with finite support, by Lemma~\ref{lem:finite-support-closed} this set is closed. Furthermore, Lemma \ref{lem:continuous} implies that the set $H_{\psi,\pi,1}$ is closed whenever $(\psi,\pi,1)$ is a nontrivial minimal valid tuple for $M_b$. The thesis now follows as $Q_b$ can be defined as set $(d)$ in Lemma \ref{lem:equal-sets}.
\endproof

For any subsets $R, P \subseteq \R^n$, define
$$V_{R,P} := \big\{(s,y) \in \R^{(\R^n)} \times \R^{(\R^n)}: s(r) = 0 \:\forall r \not\in R, \; y(p) = 0 \:\forall p \not\in P\big\}.$$
When convenient, we will see $V_{R,P}$ as a subset of $\R^R\times\R^P$ by dropping the variables set to 0.

\begin{lemma}\label{lem:closedness-canonical-face}
For any $R, P \subseteq \R^n$, $V_{R,P}$ is a closed subspace of $\R^{(\R^n)} \times \R^{(\R^n)}$ under the topology induced by $|(\cdot, \cdot)|_*$.
\end{lemma}

\proof{Proof.}
$V_{R,P}$ can be seen as the intersection of a family of hyperplanes with finite support and therefore, by Lemma~\ref{lem:finite-support-closed} this set is closed.
\endproof

Define $\cl(\cdot)$ as the closure operator with respect to the topology induced by $|(\cdot, \cdot)|_*$.

\begin{theorem}\label{thm:closure-char}
$Q_b=\cl(\conv(M_b))=\conv(M_b) + (\R^{(\R^n)}_+ \times \R^{(\R^n)}_+)$.
\end{theorem}

\proof{Proof.}
We first show that $Q_b\supseteq\cl(\conv(M_b))$.
Since, under the topology induced by $|(\cdot, \cdot)|_*$, $Q_b$ is a closed convex set by Lemma \ref{lem:Qb-closed}, it suffices to show that $Q_b\supseteq M_b$. This follows from the fact that $M_b\subseteq\R^{(\R^n)}_+ \times \R^{(\R^n)}_+$ and every inequality that defines $Q_b$ is valid for $M_b$.\smallskip

We next show that $Q_b\subseteq\cl(\conv(M_b))$.
Consider a point $(s,y) \not\in \cl(\conv(M_b))$. By the Hahn-Banach theorem, there exists a continuous linear functional that separates $(s,y)$ from $\cl(\conv(M_b))$. In other words, there exist two functions $\psi, \pi:\R^n\to \R$ and a real number $\alpha$ such that $F_{\psi,\pi}(s,y) < \alpha$ and $\cl(\conv(M_b)) \subseteq H_{\psi,\pi,\alpha}$, implying that $(\psi,\pi,\alpha)$ is a valid tuple for $M_b$. Thus, $(s,y)\notin Q_b$.\smallskip

We now show that $\conv(M_b) + (\R^{(\R^n)}_+ \times \R^{(\R^n)}_+)\subseteq Q_b$. Consider any point $(s_1, y_1) + (s_2, y_2)$, where $(s_1, y_1) \in \conv(M_b)$ and $s_2 \geq 0, y_2 \geq 0$. Since $Q_b$ can be written as the set $(d)$ in Lemma \ref{lem:equal-sets} and $\conv(M_b)\subseteq\R^{(\R^n)}_+ \times \R^{(\R^n)}_+$, we just need to verify that $(s_1, y_1) + (s_2, y_2) \in H_{\psi,\pi,1}$ for all valid tuples $(\psi, \pi,1)$ such that $\psi, \pi \geq 0$. This follows because $(s_1, y_1) \in H_{\psi, \pi,1}$ and $(s_2, y_2)$ and $\psi, \pi$ are all nonnegative.\smallskip



We finally show that $\conv(M_b) + (\R^{(\R^n)}_+ \times \R^{(\R^n)}_+)\supseteq Q_b$.
Consider $(s^*,y^*) \not\in \conv(M_b) + (\R^{(\R^n)}_+ \times \R^{(\R^n)}_+)$. We prove that $(s^*,y^*)\not\in Q_b$.
This is obvious when $(s^*,y^*)\notin\R^{(\R^n)}_+ \times \R^{(\R^n)}_+$. Therefore we assume $s^*\ge0$, $y^*\ge0$.
Let $R \subseteq \R^n$ be a finite set containing the support of $s^*$ and satisfying $\cone(R)=\R^n$ (where $\cone(R)$ denotes the conical hull of $R$),  and let $P\subseteq \R^n$ be a finite set containing the support of $y^*$. 
Then $(s^*,y^*) \not\in \conv(M_b \cap V_{R,P}) + (\R^R_+ \times \R^P_+)$. (We use the same notation $(s^*,y^*)$ to indicate the restriction of $(s^*,y^*)$ to $\R^R\times\R^P$.)
Since $M_b \cap V_{R,P}$ is the inverse image of the closed set $b + \Z^n$ under the linear transformation given by the matrix $(R,P)$, $M_b \cap V_{R,P}$ is closed in the usual finite dimensional topology of $V_{R,P}$. Therefore, by Lemma \ref{lemma:closed}, $\conv(M_b \cap V_{R,P}) + (\R^R_+ \times \R^P_+)$ is closed as well. This implies that there exists a valid inequality in $\R^R \times \R^P$ separating $(s^*,y^*)$ from $\conv(M_b \cap V_{R,P}) + (\R^R_+ \times \R^P_+)$. Since the recession cone of $\conv(M_b \cap V_{R,P}) + (\R^R_+ \times \R^P_+)$ contains $(\R^R_+ \times \R^P_+)$ and because $s^*,y^*\ge0$, this valid inequality is of the form $\sum_{r\in R}h(r)s(r) + \sum_{p\in P}d(p)y(p) \geq 1$, where $h(r) \geq 0$ for $r\in R$ and $d(p) \geq 0$ for $p\in P$.

 Now define the functions
$$\textstyle \psi(r) := \inf\left\{\sum_{r'\in R}h(r')s(r')  :  r = \sum_{r'\in R}r's(r'), \, s:R \to \R_+\right\},$$
\begin{multline*}
\textstyle \pi(p) := \inf\big\{\sum_{r'\in R}h(r')s(r') + \sum_{p'\in P}d(p')y(p') : \\
\textstyle p = \sum_{r'\in R}r's(r') + \sum_{p'\in P}p'y(p'), \, s: R \to \R_+,\; y: P \to \Z_+\big\}.
\end{multline*}
Since $\cone(R) = \R^n$, $\psi$ and $\pi$ are well-defined functions. As the sum only involves nonnegative terms, $\psi, \pi \geq 0$. It was shown in~\cite[Theorem 5]{basu2016choose} that $(\psi, \pi, 1)$ is a valid tuple for $M_b$, and since $(s^*,y^*)\notin H_{\psi,\pi,1}$,  we have  $(s^*,y^*) \not\in Q_b$.
\endproof

\subsection{The closure of $\conv(I_b)$.}

In the following, we see $\R^{(\R^n)}$ as a topological vector subspace of the space $\R^{(\R^n)}\times\R^{(\R^n)}$ endowed with the topology induced by the norm $|(\cdot,\cdot)|_*$. With a slight abuse of notation, for any $y \in \R^{(\R^n)}$ we write $|y|_* := |y(0)| + \sum_{p\in \R^n}\|p\||y(p)|$.
Also, given $\pi:\R^n\to\R$ and $\alpha\in\R$, we let $H_{\pi,\alpha} := \big\{y\in\R^{(\R^n)}: \sum_{p\in \R^n} \pi(p)y(p)\ge \alpha\big\}$.


We define $G_b := \{ y \in \R^{(\R^n)} : (0,y) \in Q_b\}.$
Since $Q_b$ can be written as the set $(c)$ in Lemma \ref{lem:equal-sets}, by Corollary \ref{cor:dominated_minimal} and Remark~\ref{rem:sup-lipschitz}, we have that
\begin{equation}\label{eq:Gb}
G_b=\R^{(\R^n)}_+ \cap \bigcap\{H_{\pi,\alpha} : (\pi,\alpha) \mbox{ minimal nontrivial liftable tuple}\}.
\end{equation}

Similar to the mixed-integer case, $\conv(I_b)\subsetneq G_b$ (this will be shown in Remark \ref{rem:Gb-strict}).

\begin{theorem}\label{thm:closure-char-Ib}
$G_b=\cl(\conv(I_b))=\conv(I_b) + \R^{(\R^n)}_+$.
\end{theorem}

\proof{Proof.}
By Theorem \ref{thm:closure-char}, $Q_b=\conv(M_b)+(\R^{(\R^n)}_+\times\R^{(\R^n)}_+)$. Since the inequality $s\ge0$ is valid for $Q_b$, by taking the intersection of $Q_b$ with the subspace $\{(s,y):s=0\}$ we obtain the equality $G_b=\conv(I_b)+\R^{(\R^n)}_+$. Furthermore, since $G_b$ coincides with the intersection of the closed set $Q_b$ with the closed subspace defined by $s = 0$ (this subspace is closed by Lemma~\ref{lem:closedness-canonical-face}), $G_b$ is a closed set. Therefore, $\cl(\conv(I_b))\subseteq G_b$.

It remains to show that $G_b\subseteq\cl(\conv(I_b))$. Consider $\bar y \not\in \cl(\conv(I_b))$. By the Hahn-Banach theorem, there exists a continuous linear functional that separates $\bar y$ from $\cl(\conv(I_b))$. In other words, there exists a function $\pi:\R^n\to \R$ and $\alpha \in \R$ such that $\sum_{r\in \R^n}\pi(r)\bar y(r) < \alpha$ and $\cl(\conv(I_b)) \subseteq \{y: \sum_{r\in \R^n}\pi(r)y(r) \geq \alpha\}$, implying that $(\pi,\alpha)$ is a valid tuple for $I_b$. We may assume without loss of generality that $(\pi,\alpha)$ is a nontrivial minimal valid tuple for $I_b$. Moreover, by Lemma~\ref{lem:cts-dir-der}, $\sup_{\epsilon > 0} \frac{\pi(\epsilon r)}{\epsilon} < \infty$ for all $r \in \R^n\setminus \{0\}$. By Corollary~\ref{cor:dominated_minimal}, $(\pi, \alpha)$ is a liftable valid tuple. Thus, $\bar y\notin G_b$ by~\eqref{eq:Gb}.
\endproof

We remark that the above theorem does not seem to follow easily from Theorem \ref{thm:closure-char}, despite the similarities in the proofs.

\section{Hamel bases, affine hulls and nonnegative representation of valid tuples.}

In finite dimensional spaces, the affine hull of any subset $C$ can be equivalently described as the set of affine combinations of points in $C$ or the intersection of all hyperplanes containing $C$. Lemma \ref{lem:aff-intrinsic-extrinsic} shows that the same holds in infinite dimension.

Before stating and proving the lemma, we give a precise definition of hyperplane in infinite dimensional vector spaces.
Given a vector space $V$ over a field $\mathbb{F}$, a subset $H \subseteq V$ is said to be a {\em hyperplane} in $V$ if there exists a nonzero linear functional $F : V \to \mathbb{F}$ and a scalar $\delta \in \mathbb{F}$ such that $H = \{v \in V: F(v) = \delta\}$.

\begin{lemma}\label{lem:aff-intrinsic-extrinsic}
Let $V$ be a vector space over a field $\mathbb F$. For every $C\subseteq V$, the set of affine combinations of points in $C$ is equal to the intersection of all hyperplanes containing $C$.
\end{lemma}

\proof{Proof.}
By possibly translating $C$, we assume without loss of generality that the set of all affine combinations of points in $C$, which we denote by $L$, 
is a linear subspace. If $x\in C$, then $x$ belongs to every hyperplane containing $C$, and therefore $L$ is contained in the intersection of all hyperplanes containing $C$.

For the reverse inclusion, let $\bar x$ be a point not in $L$. By the axiom of choice, there exists a basis $B$ of $V$ containing $\bar x$ such that $B\cap L$ is a basis of $L$. Let $F$ be the linear functional that takes value 1 on $\bar x$ and 0 on every element in $B\setminus\{\bar x\}$. Then $L\subseteq\{x:F(x)=0\}$, but $F(\bar x)=1$.
\endproof

Let $V$ be a vector space over a field $\mathbb F$. In the remainder, for any $C\subseteq V$, we will use the notation $\aff(C)$ to denote the set of affine combinations of points in the set $C$, which by the above result is equal to the intersection of all hyperplanes that contain $C$. The next proposition shows that there is no hyperplane containing $M_b$.

\begin{prop}\label{prop:aff-hull-M_b}
$\aff(M_b)=\R^{(\R^n)}\times\R^{(\R^n)}$.
\end{prop}

\proof{Proof.}
Assume by contradiction that $\aff(M_b)\subsetneq\R^{(\R^n)}\times\R^{(\R^n)}$. By Lemma~\ref{lem:aff-intrinsic-extrinsic}, there exists an equation $\sum_{r\in\R^n}\gamma(r)s(r)+\sum_{p\in\R^n}\theta(p)y(p)=\alpha$ satisfied by all points in $M_b$, where $(\gamma,\theta,\alpha)\ne(0,0,0)$.
As $\R^{(\R^n)}_+\times\R^{(\R^n)}_+$ is not contained in any hyperplane, either the valid tuple $(\gamma,\theta,\alpha)$ or the valid tuple $(-\gamma,-\theta,-\alpha)$ is nontrivial. Without loss of generality, we assume that $(\gamma,\theta,\alpha)$ is nontrivial.
Let $(\gamma',\theta',\alpha)$ be a minimal valid tuple with $\gamma'\le\gamma$ and $\theta'\le\theta$.
Note that $(\gamma',\theta')\ne(0,0)$, as $(\gamma',\theta',\alpha)$ is nontrivial.
Since $(\gamma',\theta',\alpha)$ is minimal and nontrivial, Theorem \ref{thm:non-dominated} implies that $\gamma'$ and $\theta'$ are continuous nonnegative functions. Therefore, as $(\gamma',\theta')\ne(0,0)$, there exists $\bar r\in\Q^n$ such that $\gamma'(\bar r)>0$ or $\theta'(\bar r)>0$. Assume $\gamma'(\bar r)>0$ (the other case is similar) and let $(\bar s,\bar y)\in M_b$. Then there exists a large enough integer $k>0$ such that the point $(s',\bar y)$ defined by $s'({\bar r})=\bar s({\bar r})+k$ and $s'(r)=\bar s(r)$ for $r\ne\bar r$ is in $M_b$, and \[\sum_{r\in\R^n}\gamma(r)s'(r)+\sum_{p\in\R^n}\theta(p)\bar y(p)\ge
 \sum_{r\in\R^n}\gamma'(r)s'(r)+\sum_{p\in\R^n}\theta'(p)\bar y(p)>\alpha,\]
contradicting the assumption that $\sum_{r\in\R^n}\gamma(r)s(r)+\sum_{p\in\R^n}\theta(p)y(p)=\alpha$ for all $(s,y)\in M_b$.
\endproof




The characterization of $\aff(I_b)$ is more involved and requires some preliminary notions.

\subsection{Hamel bases and the solutions to the Cauchy functional equation.}

A function $\theta:\R^n\to\R$ is {\em additive} if it satisfies the following \emph{Cauchy functional equation in $\R^n$}:
\begin{equation}\label{eq:additive}
\theta(u+v)=\theta(u)+\theta(v)\mbox{ for all }u,v\in\R^n.
\end{equation}
Note that if $\theta$ is an additive function, then
$$\theta(qx)=q\theta(x)\mbox{ for every $x\in\R^n$ and $q\in\Q$}.$$

Equation \eqref{eq:additive} has been extensively studied, see e.g.~\cite{aczel1989functional}. We summarize here the main results that we will employ.

Given any $c\in\R^n$, the linear function $\theta(x)=c^Tx$ is obviously a solution to the equation. However, these are not the only solutions. Below we describe all solutions to the equation.

A {\em Hamel basis} for $\R^n$ is a basis of the vector space $\R^n$ over the field $\Q$. In other words, a Hamel basis is a subset $B\subseteq\R^n$ such that, for every $x\in\R^n$, there exists a unique choice of a  finite subset $\{ a_1,\dots, a_t\}\subseteq B$ (where $t$ depends on $x$) and nonzero rational numbers $\lambda_1,\dots,\lambda_t$ such that
\begin{equation}\label{eq:combination}
x=\sum_{i=1}^t\lambda_i a_i.
\end{equation}
The existence of $B$ is guaranteed under the axiom of choice.

For every $ a\in B$, let $c( a)$ be a real number. Define $\theta$ as follows: for every $x\in\R^n$, if \eqref{eq:combination} is the unique decomposition of $x$, set
\begin{equation}\label{eq:theta}
\theta(x)=\sum_{i=1}^t \lambda_ic( a_i).
\end{equation}
It is easy to check that a function of this type is additive. The following theorem proves that all additive functions are of this form~\cite[Theorem 10]{aczel1989functional}.

\begin{theorem}
Let $B$ be a Hamel basis of $\R^n$. Then every additive function is of the form \eqref{eq:theta} for some choice of real numbers $c( a),\, a\in B$.
\end{theorem}

\subsection{The affine hull of $I_b$.}\label{sec:aff-hull-Ib}

The following result is an immediate extension of a result of Basu, Hildebrand and K\"oppe (see \cite[Propositions 2.2--2.3]{basu2016light}).

\begin{prop}\label{prop:affine-hull}
The affine hull of $I_b$ is described by the equations
\begin{equation}\label{eq:aff-Ib}
\sum_{p\in\R^n}\theta(p)y(p)=\theta(b)
\end{equation}
for all additive functions $\theta:\R^n\to\R$ such that $\theta(p)=0$ for every $p\in\Q^n$.
\end{prop}

\proof{Proof.} By Lemma~\ref{lem:aff-intrinsic-extrinsic}, the affine hull of $I_b$ is the intersection of all hyperplanes in $\R^{(\R^n)}$ containing $I_b$.

We first show that any equation of the form~\eqref{eq:aff-Ib} gives a hyperplane that contains $I_b$. If $y\in I_b$, then there exists $z\in\Z^n$ such that $\sum_{p\in\R^n}py(p)=b+z$. This implies that
\[\sum_{p\in\R^n}\theta(p)y(p)=\theta\Bigg(\sum_{p\in\R^n}py(p)\Bigg)=\theta(b+z)=\theta(b),\]
where the first equation comes from the additivity of $\theta$ and the integrality of $y(p)$, and the last equation from $\theta(z)=0$.
This shows that every equation of the form \eqref{eq:aff-Ib} is valid for $I_b$.

Next, we prove that any hyperplane in $\R^{(\R^n)}$ containing $I_b$ has the form~\eqref{eq:aff-Ib}. Let $\sum_{p\in\R^n}\theta(p)y(p)=\alpha$ be a hyperplane containing $I_b$. We show that $\theta$ is an additive function.
Given $p\in\R^n$, let $e_p$ denote the function such that $e_p(p)=1$ and $e_p(p')=0$ for $p'\ne p$.
Given $p_1,p_2\in\R^n$, define $y_1:=e_{p_1+p_2}+e_{b-p_1-p_2}$ and $y_2:=e_{p_1}+e_{p_2}+e_{b-p_1-p_2}$.
Since $y_1,y_2\in I_b$, $\alpha=\sum_{p\in\R^n}\theta(p)y_1(p)=\sum_{p\in\R^n}\theta(p)y_2(p)$. This shows that $\theta(p_1+p_2)=\theta(p_1)+\theta(p_2)$.
Therefore $\theta$ is additive.

Since $(\theta,\alpha)$ and $(-\theta,-\alpha)$ are valid tuples, and valid tuples are nonnegative on the rationals (see Subsection 2.1.2 in \cite{basu2016light}), it follows that $\theta(p)=0$ for every $p\in\Q^n$. Finally, since $e_b\in I_b$, we have that $\alpha=\theta(b)$.
\endproof

\begin{remark}\label{rem:Gb-strict}
Since, by the above proposition, $\conv(I_b)$ is contained in some hyperplane, $\conv(I_b)\subsetneq\conv(I_b)+\R^{(\R^n)}_+=G_b$, where the equality follows from Theorem \ref{thm:closure-char-Ib}.
\end{remark}

In the following, $e_1,\dots,e_n$ denote the vectors of the standard basis of $\R^n$. For any subset $P \subseteq \R^n$, we will use the notation $$V_P := \{y \in \R^{(\R^n)}: y(p) = 0 \;\forall p \not\in P \}.$$ When convenient, we will see $V_P$ as a subset of $\R^P$ by dropping the variables set to $0$.

\begin{prop}\label{prop:affine-hull-rational}
Let $P$ be a finite subset of $\R^n$. Then $\aff(I_b)\cap V_P$ is a rational affine subspace of $\R^P$, i.e., there exist a natural number $m \leq |P|$, a rational matrix $\Theta \in \Q^{m\times |P|}$ and a vector $d \in \R^m$ such that $\aff(I_b)\cap V_P = \{y \in \R^P: \Theta y = d\}$. Moreover, $\aff(I_b)\cap V_P = V_P$ if and only if $P\subseteq \Q^n$.
\end{prop}

\proof{Proof.}
Let $I=\{p_1,\dots,p_k\}$ be a maximal subset of vectors in $P$ such that $I\cup\{e_1,\dots,e_n\}$ is linearly independent over $\Q$, and let $B$ be a Hamel basis of $\R^n$ containing $I\cup\{e_1,\dots,e_n\}$. Note that $I=\emptyset$ if and only if $P\subseteq\Q^n$.

For every $i=1,\dots,k$, let $\theta_i$ be the additive function defined by $\theta_i(p_i)=1$ and $\theta_i(p)=0$ for every $p\in B\setminus\{p_i\}$.
Note that every $\theta_i$ is an additive function that takes value 0 on the rationals, since $\{e_1,\dots,e_n\} \subseteq B$. Moreover, $\theta_i(p) \in \Q$ for all $p \in P$. Therefore, by Proposition~\ref{prop:affine-hull}, $\sum_{p\in P}\theta_i(p)y(p) = \theta_i(b)$ is an equation satisfied by $\aff(I_b)\cap V_P$ with rational coefficients on the left hand side.
Again by Proposition~\ref{prop:affine-hull}, in order to show that these equations suffice to describe $\aff(I_b)\cap V_P$ it suffices to show the following: for every additive function $\theta$ that takes value 0 on the rationals, there exist $\lambda_1,\dots,\lambda_k\in\R$ such that $\theta(p)=\sum_{i=1}^k\lambda_i\theta_i(p)$ for every $p\in P$.

Let $\theta$ be an additive function that takes value 0 on the rationals, and define $\lambda_i:=\theta(p_i)$ for $i=1,\dots,k$.
For every $p\in P$, there exist $\bar q\in\Q^n$ and $q_1,\dots,q_k\in\Q$ such that $p=\bar q+\sum_{i=1}^kq_ip_i$. Then, since $\theta_i$ is additive and $\theta_i(\bar q)=0$, we have $\theta_i(p)=\theta_i(\sum_{j=1}^kq_jp_j)=\sum_{j=1}^kq_j\theta_i(p_j)=q_i$ for every $i=1,\dots,k$. It follows that
\[\theta(p)=\theta\Bigg(\sum_{i=1}^kq_ip_i\Bigg)=\sum_{i=1}^kq_i\theta(p_i)=\sum_{i=1}^k\theta_i(p)\lambda_i.\]
We finally observe that in the above arguments, if $I\neq \emptyset$, then we get at least one non-trivial equation corresponding to $\theta_i$, $i \in I$. Therefore, $\aff(I_b)\cap V_P = V_P$ if and only if $I=\emptyset$, which is equivalent to $P\subseteq \Q^n$.
\endproof

\subsection{Every minimal tuple is equivalent to a unique nonnegative valid tuple.}

\begin{theorem}\label{th:nonneg}
For every minimal valid tuple $(\pi,\alpha)$ for $I_b$, there exists a unique additive function $\theta:\R^n\to\R$ such that  $\theta(p)=0$ for every $p\in\Q^n$,  $(\pi',\alpha')=(\pi-\theta,\alpha-\theta(b))$ and the valid tuple is minimal and satisfies $\pi'\ge0$.
\end{theorem}

Recall that the additive functions $\theta:\R^n\to\R$ such that  $\theta(p)=0$ for every $p\in\Q^n$ are precisely the functions that define the affine hull of $I_b$ (Proposition \ref{prop:affine-hull}). Thus, the above theorem answers Open Question 2.5 in \cite{basu2016light}.
The rest of this subsection is devoted to proving Theorem \ref{th:nonneg}.

Note that if $B$ is a Hamel basis of $\R^n$ such that $e_i\in B$ for all $i\in[n]$ and $\theta$ is an additive function as in \eqref{eq:theta}, the requirement that $\theta(p)=0$ for every $p\in\Q^n$ is equivalent to $c(e_i)=0$ for $i\in[n]$.
Therefore, in order to prove the theorem, we show that given a minimal valid tuple $(\pi,\alpha)$, there exists a unique additive function $\theta$ such that $\theta(e_i)=0$ for all $i\in[n]$ and $\pi-\theta$ is a nonnegative function. If this happens, note that $(\pi-\theta,\alpha-\theta(b))$ is still a minimal valid tuple.


\begin{lemma}\label{lem:minimal}
If $(\pi,\alpha)$ is a minimal valid tuple, then $\pi$ is subadditive, $\pi(z)=0$ for every $z\in\Z^n$, and $\pi$ is periodic modulo $\Z^n$ (i.e, $\pi(p+z)=\pi(p)$ for every $p\in\R^n$ and $z\in\Z^n$).
\end{lemma}

\proof{Proof.}
We refer to the proof of Theorem 6.22 in \cite{conforti2014integer}, which however assumes the nonnegativity of $\pi$.
It is easy to check that the proof of subadditivity in \cite{conforti2014integer} does not require the nonnegativity of $\pi$.
On the contrary, the proof in \cite{conforti2014integer} that $\pi(z)=0$ for every $z\in\Z^n$ uses nonnegativity of $\pi$. However, one observes that $\pi$ must be nonnegative on the rationals (and thus on the integers), and this suffices to apply the same proof as in \cite{conforti2014integer}.
Periodicity now follows as in \cite{conforti2014integer}.
\endproof

\paragraph{Some useful results from \cite{yildiz2015integer}}

We will need some results of Y\i ld\i z and Cornu\'ejols \cite{yildiz2015integer}, which need to be slightly generalized, as only valid tuples with $\alpha=1$ are considered in \cite{yildiz2015integer}.

Let $(\pi,1)$ be a minimal valid tuple.
By Lemma 12 in \cite{yildiz2015integer} (with $f=-b$ and $S=\Z^n$), $\pi$ satisfies the generalized symmetry condition (equation (4) in \cite{yildiz2015integer}), which, by periodicity of $\pi$ modulo $\Z^n$, reads as follows:
\begin{equation}\label{eq:gen-symm}
\pi(p)=\sup_{k\in\Z_{>0}}\left\{\frac{1-\pi(b-kp)}k\right\}\quad\mbox{for all $p\in\R^n$}.
\end{equation}
Then, by Proposition 17 in \cite{yildiz2015integer} (with $f=-b$, $S=\Z^n$, $X=\{0\}$), for any $p \in \R^n$,
the supremum in \eqref{eq:gen-symm} is attained if and only if $\pi(p)+\pi(b-p)=1$.
Proposition 18 in \cite{yildiz2015integer} then implies the following: if $p\in\R^n$ is such that $\pi(p)+\pi(b-p)>1$, then
\[\limsup_{k\in\Z_{>0},k\to\infty}\frac{\pi(kp)}{k}=\limsup_{k\in\Z_{>0},k\to\infty}\frac{-\pi(-kp)}{k}.\]

One straightforwardly (and patiently) verifies that when $(\pi,\alpha)$ is a minimal valid tuple with $\alpha$ not restricted to be 1, the above result generalizes as follows:
\begin{prop}\label{prop:yild-corn}
Let $(\pi,\alpha)$ be a minimal valid tuple. If $p\in\R^n$ is such that $\pi(p)+\pi(b-p)>\alpha$, then
\[\limsup_{k\in\Z_{>0},k\to\infty}\frac{\pi(kp)}{k}=\limsup_{k\in\Z_{>0},k\to\infty}\frac{-\pi(-kp)}{k}.\]
\end{prop}

\paragraph{Construction of $\theta$}

In what follows, we will assume $(\pi,\alpha)$ is a minimal valid tuple for $I_b$. Let $B$ be a Hamel basis of $\R^n$ containing the unit vectors $e_1,\dots,e_n$.
For every $ a\in B$, define
\begin{equation}\label{eq:c}
c( a):=\inf_{k\in\Z_{>0}}\frac{\pi(k a)}{k}.
\end{equation}
We will show that this is the correct choice for the constant $c( a)$.

\begin{lemma}\label{lemma:inf-sup}
For all $a\in B$, the value of $c( a)$ is finite and
\[c( a)=\inf_{k\in\Z_{>0}}\frac{\pi(k a)}{k}=\sup_{k\in\Z_{>0}}\frac{-\pi(-k a)}{k}.\]
\end{lemma}

\proof{Proof.}
We prove a sequence of claims.

\begin{claim}
The inequality ``$\inf\ge\sup$'' holds and both terms are finite.
\end{claim}

\begin{claimproof}
Let $h,k$ be positive integers. Then, by subadditivity,
$h\pi(k a)+k\pi(-h a)\ge\pi(0)=0$,
thus $\frac{\pi(k a)}{k}\ge-\frac{\pi(-h a)}{h}$. Since this holds for all positive integers $h,k$, the claim is proven.
\end{claimproof}

We now assume by contradiction that
\[\inf_{k\in\Z_{>0}}\frac{\pi(k a)}{k}-\sup_{k\in\Z_{>0}}\frac{-\pi(-k a)}{k}\ge\epsilon\]
for some $\epsilon>0$. In other words,
\begin{equation}\label{epsilon}
\inf_{k\in\Z_{>0}}\frac{\pi(k a)}{k}+\inf_{k\in\Z_{>0}}\frac{\pi(-k a)}{k}\ge\epsilon.
\end{equation}

\begin{claim}
The following equation holds:
\begin{equation}\label{inf}
\inf_{k\in\Z_{>0}}\frac{\pi(k a)}{k}+\inf_{k\in\Z_{>0}}\frac{\pi(-k a)}{k} =
\inf_{k\in\Z_{>0}}\frac{\pi(k a)+\pi(-k a)}{k}.
\end{equation}
\end{claim}

\begin{claimproof}
Since the inequality ``$\le$'' is obvious, we prove the reverse inequality.
To do so, it is sufficient to show that given positive integers $h,k$, there exists a positive integer $\ell$ such that
\begin{equation}\label{ell}
\frac{\pi(h a)}{h}+\frac{\pi(-k a)}{k}\ge\frac{\pi(\ell a)+\pi(-\ell a)}{\ell}.
\end{equation}
Choose $\ell=hk$. Then, by subadditivity,
\[k\pi(h a)+h\pi(-k a)\ge\pi(\ell a)+\pi(-\ell a).\]
After dividing by $\ell=hk$, we obtain \eqref{ell} and the claim is proven.
\end{claimproof}

By the previous claim, assumption \eqref{epsilon}  is equivalent to
\begin{equation}\label{eq:reformulate-contra}\frac{\pi(k a)}{k}+\frac{\pi(-k a)}{k}\ge\epsilon\:\mbox{ for all positive integers $k$}.\end{equation}

\begin{claim}\label{claim:strict-symmetry}
There exists a positive integer $k$ such that $\pi(k a)+\pi(b-k a)>\alpha$.
\end{claim}

\begin{claimproof}
By subadditivity, for every integer $k$ we have
\[\pi(b-k a)\ge\pi(-k a)-\pi(-b).\]
Combined with~\eqref{eq:reformulate-contra}, we get that
\[\pi(ka) + \pi(b-ka) \geq \pi(ka) +\pi(-ka) -\pi(-b) \geq \epsilon k -\pi(-b),\] for all positive integers $k$.
The right-hand side is greater than $\alpha$ if $k>\frac{\alpha+\pi(-b)}{\epsilon}$. 
\end{claimproof}


Define $p:=\bar k a$, where $\bar k$ is the positive integer guaranteed by Claim~\ref{claim:strict-symmetry}. 

\begin{claim}\label{claim:final}The following inequalities hold:
\[\limsup_{k\in\Z_{>0},k\to\infty}\frac{\pi(kp)}{k}\ge\bar k\cdot\inf_{k\in\Z_{>0}}\frac{\pi(k a)}{k},\quad
\limsup_{k\in\Z_{>0},k\to\infty}\frac{-\pi(-kp)}{k}\le\bar k\cdot\sup_{k\in\Z_{>0}}\frac{-\pi(-k a)}{k}\]
\end{claim}

\begin{claimproof}
Since the $\limsup$ is always at least as large as the inf,
\[\limsup_{k\in\Z_{>0},k\to\infty}\frac{\pi(kp)}{k}\ge\inf_{k\in\Z_{>0}}\frac{\pi(kp)}{k}
=\bar k\cdot\inf_{k\in\Z_{>0}}\frac{\pi(k\bar k a)}{k\bar k}
\ge\bar k\cdot\inf_{h\in\Z_{>0}}\frac{\pi(h a)}{h}\]
and thus the first inequality is verified.

Since the $\limsup$ is always at most as large as the sup,
\[\limsup_{k\in\Z_{>0},k\to\infty}\frac{-\pi(-kp)}{k}\le\sup_{k\in\Z_{>0}}\frac{-\pi(-kp)}{k}
=\bar k\cdot\sup_{k\in\Z_{>0}}\frac{-\pi(-k\bar k a)}{k\bar k}
\le\bar k\cdot\sup_{h\in\Z_{>0}} \frac{-\pi(-h a)}{h}\]
and thus the second inequality holds.
\end{claimproof}

Since $\pi(p)+\pi(b-p)>\alpha$, by Proposition \ref{prop:yild-corn},
\[\limsup_{k\in\Z_{>0},k\to\infty}\frac{\pi(kp)}{k}=\limsup_{k\in\Z_{>0},k\to\infty}\frac{-\pi(-kp)}{k}.\]
Claim~\ref{claim:final} then implies that $\inf_{k\in\Z_{>0}}\frac{\pi(k a)}{k}\le \sup_{k\in\Z_{>0}}\frac{-\pi(-k a)}{k}$. But this contradicts~\eqref{eq:reformulate-contra}. This concludes the proof of the lemma.
\endproof

Now let $\theta$ be defined as in \eqref{eq:theta}, where the constants $c( a)$ for $ a\in B$ are chosen as in \eqref{eq:c}. In the next two lemmas we prove that $\theta(e_i)=0$ for all $i\in[n]$ and $\pi-\theta$ is nonnegative.

\begin{lemma}
$\theta(e_i)=0$ for all $i\in[n]$.
\end{lemma}

\proof{Proof.}
Fix $i\in[n]$. Since $e_i\in B$, it is sufficient to check that $c(e_i)=0$. By \eqref{eq:c}, $c(e_i)=\inf_{k\in\Z_{>0}}\frac{\pi(ke_i)}{k}$.
Since $\pi(ke_i)=0$ for all $k\in \Z$ by Lemma \ref{lem:minimal}, we immediately see that $c(e_i)=0$.
\endproof

\begin{lemma}
If $\theta$ is defined as in \eqref{eq:theta}, with the constants $c(a)$ given in \eqref{eq:c}, then the function $\pi-\theta$ is nonnegative.
\end{lemma}

\proof{Proof.}
Let $x\in\R^n$. Then there exist $ a_1,\dots, a_t\in B$ and nonzero rational numbers $\lambda_1,\dots,\lambda_t$ such that
$x=\sum_{i=1}^t\lambda_i a_i$, and we have $\theta(x)=\sum_{i=1}^t\lambda_ic( a_i)$. We prove that $\pi(x)-\theta(x)\ge0$.

For every $i\in\{1,\dots,t\}$, we can write $\lambda_i=\frac{p_i}{q_i}$, where every $p_i$ is a nonzero integer and every $q_i$ is a positive integer. Define $Q:=\prod_{j=1}^t q_j$. Take arbitrary positive integers $k_1,\dots,k_t$ (these numbers will be fixed later) and define $K:=\prod_{j=1}^tk_j$. Since $\frac{Q}{q_i}$ and $\frac{K}{k_i}$ are positive integers for every $i$, by subadditivity we have
\[QK\pi(x)+\sum_{i=1}^t\frac{QK}{q_ik_i}\pi(-k_ip_i a_i) \ge
\pi\bigg(QKx-\sum_{i=1}^t QK\lambda_i a_i\bigg)=\pi(0)=0.\]
This implies that
\begin{equation}\label{eq:nonneg}
\pi(x)\ge\sum_{i=1}^t\lambda_i\,\frac{-\pi(-k_ip_i a_i)}{k_ip_i}.
\end{equation}

Now fix $\epsilon>0$.
If $i$ is an index such that $p_i>0$, by Lemma \ref{lemma:inf-sup} we can choose $k_i$ such that $\frac{-\pi(-k_i a_i)}{k_i}\ge c( a_i)-\epsilon$. Then by subadditivity
\[\frac{-\pi(-k_ip_i a_i)}{k_ip_i}\ge\frac{-\pi(-k_i a_i)}{k_i}\ge c( a_i)-\epsilon.\]
If $i$ is an index such that $p_i<0$, by Lemma \ref{lemma:inf-sup} we can choose $k_i$ such that $\frac{\pi(k_i a_i)}{k_i}\le c( a_i)+\epsilon$. Then by subadditivity
\[\frac{-\pi(-k_ip_i a_i)}{k_ip_i}\le\frac{\pi(k_i a_i)}{k_i}\le c( a_i)+\epsilon.\]
Then, rembering that $\lambda_i>0$ if and only if $p_i>0$, equation \eqref{eq:nonneg} gives
$\pi(x)\ge \sum_{i=1}^t \lambda_i c( a_i)-\epsilon(\sum_{i=1}^t |\lambda_i| )$.
Since this holds for every $\epsilon>0$, we have $\pi(x)\ge \sum_{i=1}^t \lambda_ic( a_i)$ and thus $\pi(x)-\theta(x)\ge0$.
\endproof

This concludes the proof of the existence of $\theta$. To conclude, it only remains to show that the choice of $\theta$ is unique.
To see this, let $\theta'$ be any additive function such that $\pi-\theta'$ is nonnegative.
For every $a\in B$ and $k\in \Z_{>0}$, we have
\[\pi(ka)\ge\theta'(ka)=k\theta'(a),\qquad \pi(-ka)\ge\theta'(-ka)=-k\theta'(a).\]
This implies
\[\sup_{k\in\Z_{>0}}\frac{-\pi(-ka)}k\le\theta'(a)\le\inf_{k\in\Z_{>0}}\frac{\pi(ka)}k.\]
Lemma \ref{lemma:inf-sup} then shows that $\theta'(a)=c(a)=\theta(a)$, which proves the uniqueness of $\theta$.


\section{Canonical faces and recession cones.}\label{sec:rec-cones-faces}



A \emph{canonical face} of $\conv(M_b)$ is a face of the form $F=\conv(M_b)\cap V_{R,P}$ for some $R,P\subseteq\R^n$. 
If $R$ and $P$ are finite, $F$ is a \emph{finite canonical face} of $\conv(M_b)$. The same definitions can be given for $\conv(I_b)$, $\cl(\conv(M_b))$ and $\cl(\conv(I_b))$. It is not hard to see that $\conv(M_b) \cap V_{R,P} = \conv(M_b \cap V_{R,P})$ and $\conv(I_b) \cap V_{P} = \conv(I_b \cap V_{P})$.
This observation shows that the finite canonical faces of $\conv(M_b)$ and $\conv(I_b)$ are finite dimensional integer programming models.


The notion of recession cone of a closed convex set is standard (see, e.g., \cite{lemarechal1996convex}). We extend it to general convex sets in general (possibly infinite-dimensional) vector spaces in the following way. Let $V$ be a vector space and let $C\subseteq V$ be a convex set. For any $x \in C$, define $$C_\infty(x) := \{r \in V: x + \lambda r \in C \textrm{ for all } \lambda \geq 0\}.\footnote{Using the Hahn-Banach separation theorem, it can be shown that if $V$ is a topological vector space and $C$ is a closed convex subset, then $C_\infty(x) = C_\infty(x')$ for all $x, x' \in C$.}$$
We define the \emph{recession cone} of $C$ as $\rec(C) := \bigcap_{x\in C}C_\infty(x).$ Theorem \ref{thm:closure-char} yields the following result.

\begin{cor}\label{cor:canonical} Given $R, P\subseteq \R^n$ (not necessarily finite), we have that $ \conv(M_b)\cap V_{R,P}=\cl(\conv(M_b)) \cap V_{R,P}$ if and only if $\rec(\conv(M_b)\cap V_{R,P})=(\R^{(\R^n)}_+\times \R^{(\R^n)}_+) \cap V_{R,P}$.

Given $P\subseteq \R^n$, we have that $ \conv(I_b)\cap V_{P}=\cl(\conv(I_b)) \cap V_{P}$ if and only if  $ \rec(\conv(I_b)\cap V_{P})=\R^{(\R^n)}_+ \cap V_P$.


\end{cor}


\proof{Proof.}
Let $F=\conv(M_b)\cap V_{R,P}$ be a canonical face of $\conv(M_b)$. By Theorem~\ref{thm:closure-char},
\begin{align*}
\cl(\conv(M_b)) \cap V_{R,P} & =  \big(\conv(M_b) + (\R^{(\R^n)}_+ \times \R^{(\R^n)}_+)\big) \cap V_{R,P} \\
& =  (\conv(M_b) \cap V_{R,P}) + ((\R^{(\R^n)}_+\times \R^{(\R^n)}_+)\cap V_{R,P}) \\
& =  F + ((\R^{(\R^n)}_+\times \R^{(\R^n)}_+)\cap V_{R,P}) .
\end{align*}
Then the result  follows from the fact that $F = F + ((\R^{(\R^n)}_+\times \R^{(\R^n)}_+)\cap V_{R,P})$ if and only if the recession cone of $F$ is $((\R^{(\R^n)}_+\times \R^{(\R^n)}_+)\cap V_{R,P})$.

The proof for $I_b$ is the same, where we use Theorem~\ref{thm:closure-char-Ib} instead of Theorem~\ref{thm:closure-char}.
\endproof

Given  $P \subseteq \R^n$, let $C^P: = \conv(I_b) \cap V_P$ be  the canonical face of $\conv(I_b)$ corresponding to $P$. Also define $L$ to be the linear space parallel to $\aff(\conv(I_b))$. Observe that, by  Proposition~\ref{prop:affine-hull}, $L$ is the set of all $y \in \R^{(\R^n)}$ that satisfy $\sum_{p\in \R^n}\theta(p)y(p) = 0$ for all additive functions $\theta: \R^n \to \R$ such that $\theta(p)=0$ for all $p\in\Q^n$.

\begin{theorem}\label{thm:rec-cone-finite}
For every finite subset $P\subseteq\R^n$, the following are all true:
\begin{enumerate}[(a)]
\item the finite canonical face $C^P$ is a rational polyhedron in $\R^P$;\item every extreme ray of $C^P$ is spanned by some $d \in \Z^P_+$ such that $\sum_{p\in P}pd(p)\in\Z^n;$
\item $\rec(C^P) = L \cap \R^{(\R^n)}_+\cap V_P = (L\cap V_P)\cap \R^P_+$ whenever $C^P\neq \emptyset$;
\end{enumerate}
\end{theorem}

\proof{Proof.}
By dropping variables set to zero,
$I_b\cap V_P$ is the set of vectors $y\in\Z^P_+$ such that
$\sum_{p\in P}py(p)\in b+\Z^n.$
We say that a feasible point $y\in I_b\cap V_P$ is minimal if there is no feasible point $y'\ne y$ such that $y'\le y$.
Every vector $d\in\Z^P_+\setminus \{0\}$  such that
$\sum_{p\in P}pd(p)\in\Z^n$ is called a {\em ray}. A ray $d$ is minimal if there is no ray $d'\ne d$ such that $d'\le d$.

Given feasible points $y$, $y'$ such that  $y'\ne y$, $y'\le y$, we have that $y-y'$ is a ray. This shows that any feasible point $y$ is the sum of a minimal feasible point and a ray. Further, given feasible rays $d$, $d'$ such that  $d'\ne d$, $d'\le d$, we have that $d-d'$ is also a ray. Thus, given a ray $d$, we can express it as the sum of a minimal ray $d'$ and a ray $d-d'$ whose coordinates are smaller than $d$. Iterating a finite number of times (since we have nonnegative integer coordinates), we can express $d$ as a nonnegative integer combination of minimal rays. Therefore, every feasible point $y$ is the sum of a minimal feasible point $\bar y$ and a nonnegative integer combination of minimal rays.

Let $Y$ be the set of points that are the sum of a minimal feasible point and a nonnegative integer combination of minimal rays.
The above observation proves that $Y=I_b\cap V_P$. By Gordan--Dickson lemma (see, e.g., \cite{dickson1913finiteness}), the set of minimal feasible points and the set of minimal rays are both finite, i.e. there exist finite sets $E\subseteq \Z^P_+$ and $R \subseteq \Z^P_+$ such that $Y = E + \integ.\cone(R)$, where $\integ.\cone(R)$ denotes the set of all nonnegative integer combinations of vectors in $R$. Since $\conv(I_b)\cap V_P=\conv(Y) = \conv(E + \integ.\cone(R)) = \conv(E) + \conv(\integ.\cone(R)) = \conv(E) + \cone(R)$, where $\cone(R)$ denotes the conical hull of $R$,  by the Minkowski-Weyl Theorem~\cite[Theorem 3.13]{conforti2014integer} we have that $\conv(I_b)\cap V_P$ is a rational polyhedron.


The above analysis proves $(a)$ and $(b)$ simultaneously. We now prove $(c)$.

We first show that $\rec(C^P) \subseteq L \cap \R^{(\R^n)}_+\cap V_P.$ Consider any $\bar d \in \rec(C^P)$. By part $(b)$, $\bar d$ is a nonnegative combination of vectors $d\in \Z^P_+$ such that $\sum_{p\in P}pd(p)\in\Z^n.$ Observe that each such $d$ is in $L$ because for any additive function $\theta$ which is 0 on the rationals, we obtain $0 = \theta(\sum_{p\in P}pd(p)) = \sum_{p\in P}\theta(p)d(p)$. Thus, $\bar d \in L$ since $L$ is a linear space. Therefore, $\rec(C^P) \subseteq L \cap \R^{(\R^n)}_+\cap V_P.$

We now want to establish that $L \cap \R^{(\R^n)}_+\cap V_P\subseteq\rec(C^P)$. First, consider any $d \in L \cap \R^{(\R^n)}_+\cap V_P$ such that $d \in \Q^P$. Let $\lambda > 0$ be such that $\bar d := \lambda d \in \Z^P_+$. We claim that $\sum_{p\in P}p\bar d(p)\in\Q^n.$ Otherwise, there exists\footnote{Such an additive function can be constructed by first constructing a Hamel basis of $\R^n$ containing $\sum_{p\in P}p\bar d(p), e_1, \ldots, e_n$, and setting $\theta$ to be $1$ on $\sum_{p\in P}p\bar d(p)$ and $0$ everywhere else on this basis.} an additive function $\theta:\R^n \to \R$ that takes value zero on the rationals such that $0 \neq\theta(\sum_{p\in P}p\bar d(p)) = \sum_{p\in P}\theta(p)\bar d(p) = \lambda \sum_{p\in P}\theta(p) d(p)$, which violates the hypothesis that $d\in L$. Since $\sum_{p\in P}p\bar d(p)\in\Q^n$, this implies that there exists a positive scaling $\tilde d$ of $d$ such that $\sum_{p\in P}p \tilde d(p)\in\Z^n$. It is easy to verify that $\tilde d \in \rec(C^P)$ whenever $C^P\neq \emptyset$, and therefore $d \in \rec(C^P)$. This shows that all rational vectors in $L \cap \R^{(\R^n)}_+\cap V_P$ are in $\rec(C^P)$. Since, by Proposition \ref{prop:affine-hull-rational}, $L\cap V_P$ is a rational subspace of $\R^P$, we conclude that $L \cap \R^{(\R^n)}_+\cap V_P\subseteq\rec(C^P)$.
\endproof

Note that the above theorem holds even if $P\not\subseteq\Q^n$. In particular, a corner polyhedron is a rational polyhedron even if $P\not\subseteq \Q^n$. The next result characterizes the corner polyhedra with $P\subseteq\Q^n$.

\begin{theorem}\label{thm:rational-P}
Let $P\subseteq\R^n$ be finite such that $C^P\ne\emptyset$. Then the following are equivalent:
\begin{enumerate}[(a)]
\item $P\subseteq\Q^n$;
\item $\rec(C^P) = \R^P_+$;
\item the dimension of $C^P$ is $|P|$;
\item $C^P=G_b\cap V_P$.
\end{enumerate}
\end{theorem}

\proof{Proof.} If we assume $(a)$ and $C^P \neq \emptyset$, we must have $b \in \Q^n$. This implies that $L = \aff(I_b)$, since $\theta(b) = 0$ for all additive functions that are 0 on the rationals. Now $(a)$ implies $(b)$ by Proposition~\ref{prop:affine-hull-rational} and Theorem~\ref{thm:rec-cone-finite}$(c)$. $(b)$ clearly implies $(c)$. $(c)$ implies $(a)$ by Proposition~\ref{prop:affine-hull-rational}.
The equivalence of $(b)$ and $(d)$ follows from Corollary \ref{cor:canonical} and Theorem~\ref{thm:closure-char-Ib}. 
\endproof

We finally characterize the recession cone of $\conv(I_b)$.

\begin{cor}\label{cor:rec-I_b} $\conv(I_b)_{\infty}(x) = L \cap \R^{(\R^n)}_+$ for every $x \in \conv(I_b)$. Consequently, $\rec(\conv(I_b)) = L \cap \R^{(\R^n)}_+$.
\end{cor}
\proof{Proof.}
%
%
%
Fix $x\in\conv(I_b)$ and $y\in\R^{(\R^n)}_+$. Let $P$ be the union of the supports of $x$ and $y$. Then $y\in\conv(I_b)_{\infty}(x)$ if and only if $y\in C^P_\infty(x)=\rec(C^P)$, which, by Theorem~\ref{thm:rec-cone-finite}$(c)$, happens if and only if $y\in  L \cap \R^{(\R^n)}_+ \cap V_P$, and this is equivalent to $y\in  L \cap \R^{(\R^n)}_+$.
\endproof

\paragraph{Examples.} We now give examples that illustrate some of the pathologies that can show up.

\begin{example}\label{ex:not-closed} There are finite dimensional faces of $\conv(M_b)$ that are not closed (in the finite dimensional topology).
Let $n=1$, $b\in\Q$, $\omega\in\R\setminus\Q$, $R=\{-1\}$, $P=\{b,\omega\}$.
Consider the point $(\bar s,\bar y)$ defined by $\bar s(-1):=0$ and $\bar y(b):=\bar y(\omega)=1$.
Note that $(\bar s,\bar y)\notin\conv(M_b)\cap V_{R,P}$, as the only point in $M_b$ satisfying $s(-1)=0$ and $y(b)\le1$ has $y(b)=1,\,y(\omega)=0$.

We now show that $(\bar s,\bar y)\in\cl(\conv(M_b)\cap V_{R,P})$ by constructing for every $\epsilon>0$ a point in $\conv(M_b)\cap V_{R,P}$ whose Euclidean distance from $(\bar s,\bar y)$ is at most $\epsilon$.
So fix $\epsilon>0$. Let $\hat y(\omega)$ be a positive integer such that the fractional part of $\omega \hat y(\omega)$ is at most $\epsilon$.
Let $\hat s(-1)$ be equal to this fractional part, and $\hat y(b)=1$. Then $(\hat s,\hat y)\in M_b\cap V_{R,P}$. Let $(\tilde s,\tilde y)$ be the point of $M_b\cap V_{R,P}$ defined by $\tilde y(b):=1,\,\tilde s(-1):=\tilde y(\omega):=0$.
Then the point $\frac{1}{\hat y(\omega)}(\hat s,\hat y)+ \frac{\hat y(\omega)-1}{\hat y(\omega)}(\tilde s,\tilde y)$ is in $\conv(M_b)\cap V_{R,P}$ and its distance from $(\bar s,\bar y)$ is $\frac{\epsilon}{\hat y(\omega)}\le\epsilon$.
\end{example}

\begin{remark}\label{rem:Q_b-strict}
Since $Q_b=\cl(\conv(M_b))$ by Theorem \ref{thm:closure-char}, for every $R,P\subseteq\R^n$ the set $Q_b\cap V_{R,P}$ is closed by Lemma \ref{lem:closedness-canonical-face}. The previous example gives sets $R,P$ such that $\conv(M_b)\cap V_{R,P}$ is not closed. Thus $\conv(M_b)$ is a strict subset of $Q_b$.
\end{remark}

\begin{example}\label{ex:pure-integer} We established in Theorem~\ref{thm:closure-char-Ib} and Remark~\ref{rem:Gb-strict} that \begin{equation}\label{eq:strict-subset-example}\conv(I_b) \subsetneq \left\{y \in \R^{(\R^n)}: \sum_{r\in \R^n} \pi(r)y(r) \geq \alpha, \;\; (\pi, \alpha) \; \textrm{ valid tuple for }I_b \right\} = \cl(\conv(I_b)).\end{equation} We now give a concrete example to illustrate this phenomenon. Let $n=1$, $b \in \Q$, $\omega \in \R\setminus \Q$ and $P = \{b, \omega, 1 - \omega\}.$ Consider the canonical face of $I_b$ given by $C^P = \conv(I_b) \cap V_P$, i.e.,
$$C^P = \conv\{(y_\omega, y_{1-\omega}, y_b) \in \Z^3_+: \omega y_\omega + (1-\omega)y_{1-\omega} + by_b \in b + \Z\}.$$
Since $\omega \not\in \Q$, it follows that for any $\bar y \in C^P$ we must have $\bar y_\omega = \bar y_{1-\omega}$. This implies that the recession cone of $C^P$ is not full-dimensional. Moreover, we know that $\rec(C^P) = (L \cap V_P) \cap \R^P_+$ by Theorem~\ref{thm:rec-cone-finite} $(c)$. In this case, considering a Hamel basis of $\R$ which includes $\omega$ and the particular additive function $\theta$ obtained by setting $\theta(\omega) = 1 = -\theta(1-\omega)$ and $0$ on the remaining elements of the Hamel basis, we get the desired finite dimensional recession cone.

Also, observe that by Corollary~\ref{cor:canonical}, $\cl(\conv(I_b))\cap V_P$ has $\R^P_+$ as its recession cone. Thus, this example explicitly shows the strict containment in~\eqref{eq:strict-subset-example}.
\end{example}

\section{Sufficiency of nontrivial minimal liftable functions to describe corner polyhedra.}\label{sec:conv-Ib-nonneg}

\begin{lemma}\label{lem:aff-rec-poly}
Given a closed, convex set $X \subseteq \R^d$, let $L_X$ be the linear space parallel to  $\aff(X)$. If $\rec(X) = L_X \cap \R^d_+$, then $(X + \R^d_+) \cap \aff(X) = X$.  
\end{lemma}
\proof{Proof.} The inclusion $\supseteq$ is clear. Now consider any $x \in X$ and $y \in \R^d_+$ such that $x + y \in \aff(X)$. Thus, $y \in \aff(X) - x = L_X$ since $x \in X$. Therefore $y \in L_X \cap \R^d_+ = \rec(X)$. Therefore, $x + y \in X$. 
\endproof

\begin{theorem}\label{thm:conv-Ib-nonneg} 
Let $C^P$ be a corner polyhedron for some finite set $P \subseteq \R^n$. Then $C^P = G_b \cap \aff(C^P)$.
\end{theorem}

\proof{Proof.} If $C^P = \emptyset$, the equality is trivial since $\aff(C^P) = \emptyset$. So we assume $C^P \neq \emptyset$. Since $C^P \subseteq V_P$, we have 
$$
\begin{array}{rcl}
G_b \cap \aff(C^P) & = & (G_b \cap V_P) \cap  \aff(C^P) \\
 & = & ((\conv(I_b) + \R_+^{(\R^n)})  \cap V_P) \cap  \aff(C^P) \\
& = & (\conv(I_b) \cap V_P + \R_+^{P}) \cap  \aff(C^P) \\
& = & (C^P + \R_+^P) \cap \aff(C^P)
\end{array}
$$
where the second equality follows Theorem \ref{thm:closure-char-Ib}, the third equality follows from the observation that $(\conv(I_b) + \R_+^{(\R^n)}) \cap V_P = (\conv(I_b) \cap V_P + \R_+^{P})$. We verify below that $\rec(C^P) = L' \cap \R^P_+$, where $L'$ is the linear space parallel to $\aff(C^P)$. Applying Lemma~\ref{lem:aff-rec-poly} with $X = C^P$, we obtain that $G_b \cap \aff(C^P) = (C^P + \R_+^P) \cap \aff(C^P) = C^P$, as desired.

We finally verify that $\rec(C^P) = L' \cap \R^P_+$. Since $C^P \neq\emptyset$ and $C^P \subseteq \R^P_+$, we have that $\rec(C^P) \subseteq L' \cap \R^P_+$. Indeed, consider any $x \in C^P$. Then, $x + \rec(C^P) \subseteq C^P \subseteq \aff(C^P)$ and so $\rec(C^P) \subseteq \aff(C^P) - x = L'$. To see the reverse inclusion, observe that $C^P=\conv(I_b)\cap V_P = \conv(I_b \cap V_P)$. Thus, $\aff(C^P) = \aff(\conv(I_b\cap V_P)) = \aff(I_b \cap V_P) \subseteq \aff(I_b) \cap V_P$. Therefore, $L' = \aff(C^P) - x \subseteq (\aff(I_b) \cap V_P) - x = L \cap V_P$, where $L$ is the linear space parallel to $\aff(I_b)$, and $x \in C^P$ is an arbitrary point. Hence, $L ' \cap \R^P_+ \subseteq (L \cap V_P) \cap \R^P_+ = \rec(C^P)$, where we use Theorem~\ref{thm:rec-cone-finite} (c). We are done establishing the reverse inclusion.
%
\endproof

\begin{cor}\label{cor:conv-Ib-nonneg} Let $\cl(\cdot)$ be the closure operator with respect to the topology induced by $|\cdot|_*$. Then
$$\conv(I_b) = G_b \cap \aff(I_b) = \cl(\conv(I_b)) \cap \aff(I_b) = (\conv(I_b)+\R^{(\R^n)}_+) \cap \aff(I_b).$$
\end{cor} 

\proof{Proof.} The first equality follows from Theorem~\ref{thm:conv-Ib-nonneg} and the facts that $$\conv(I_b) =\bigcup_{P\mbox{ \small finite subset of } \R^n}C^P, \qquad \aff(I_b) = \bigcup_{P\mbox{ \small finite subset of } \R^n}\aff(C^P)$$ because $\conv(I_b)$ is a set of finite support functions. The remaining equalities are a consequence of Theorem~\ref{thm:closure-char-Ib}.
\endproof


\begin{theorem}\label{thm:facet-restriction} Let $C^P$ be a nonempty corner polyhedron where  $P \subseteq \Q^n$. Any valid inequality for $C^P$ that is not valid for $\R^P_+$, is dominated by the restriction  to $\R^P$ of the inequality defined by some nontrivial, minimal liftable tuple for $I_b$.
\end{theorem}
\proof{Proof.} By Theorem~\ref{thm:rational-P} $(b)$, $\rec(C^P) = \R^P_+$; thus any inequality that is valid for $C^P$ but not for $\R^P_+$ is of the form $\sum_{p \in P}d(p) y(p) \geq 1$, where $d(p)$, $p \in P$ are nonnegative coefficients. Let $R = \{\pm e^1, \ldots, \pm e^n\}$, where $e^i$, $i=1, \ldots, n$ are the standard unit vectors in $\R^n$; note that $\cone(R) = \R^n$. Then, $C^P$ is a face of the canonical face $F = \conv(M_b) \cap V_{R,P}$ of $M_b$. Since $R,P$ are sets of rational vectors and $C^P\ne \emptyset$, $b$ is rational as well. Therefore by Meyer's theorem~\cite{Meyer} (see also~\cite[Theorem 4.30]{conforti2014integer})  both $F$ and $C^P$ are rational polyhedra. Since $\sum_{p \in P}d(p) y(p) \geq 1$ is a valid inequality for the face $C^P$ of $F$, it can be lifted to a valid inequality $\sum_{r \in R}h(r)s(r) + \sum_{p \in P}d(p) y(p)\geq 1$ for the polyhedron $F$. As in the proof of Theorem~\ref{thm:closure-char}, using~\cite[Theorem 5]{basu2016choose} one can find a minimal valid tuple $(\psi, \pi, 1)$ for $M_b$ that defines an inequality whose restriction  dominates $\sum_{r \in R}h(r)s(r) + \sum_{p \in P}d(p) y(p) \geq 1$. In particular, the restriction of the inequality defined by the liftable tuple $(\pi, 1)$ dominates the inequality $\sum_{p \in P}d(p) y(p) \geq 1$.
\endproof

\section*{Acknowledgement.} A. Basu and J. Paat were supported by the NSF grant CMMI1452820. 
M. Conforti and M. Di Summa were supported by the Italian national grant ``PRIN 2015'' and the grant ``SID 2016'' of the University of Padova.
The authors thank two anonymous reviewers for very insightful comments. This helped to remove errors, simplify some proofs and also strengthen some results from a previous version of the paper.

\appendix

\section{Missing proofs.}

\proof{Proof of Proposition~\ref{prop:psi_sublinear}}
Assume that $\psi$ is not subadditive. Then $\psi(r_1+r_2)>\psi(r_1)+\psi(r_2)$ for some $r_1, r_2\in \R^n$.
Let $\psi':\R^n\to \R$ be defined as $\psi'(r_1+r_2):=\psi(r_1)+\psi(r_2)$ and $\psi'(r):=\psi(r)$ for every $r\ne r_1+r_2$.
Then $(\psi',\pi,\alpha)$ is easily seen to be a valid tuple, a contradiction to the minimality of $(\psi,\pi,\alpha)$.

Now assume that $\psi$ is not positively homogenous. Then $\psi(\lambda r_1)<\lambda\psi(r_1)$ for some $r_1\in\R^n$ and $\lambda>0$.
Let $\psi':\R^n\to \R$ be defined as $\psi'(r_1):=\frac{\psi(\lambda r_1)}{\lambda}$ and $\psi'(r):=\psi(r)$ for every $r\ne r_1$.
Again, $(\psi',\pi,\alpha)$ is a valid tuple, a contradiction to the minimality of $(\psi,\pi,\alpha)$. Thus $\psi$ is sublinear.

Finally, assume that $\pi(r_1)>\psi(r_1)$ for some $r_1\in\R^n$. Let $\pi':\R^n\to \R$ be defined as $\pi'(r_1)=\psi(r_1)$ and $\pi'(r)=\pi(r)$ for every $r\ne r_1$. The tuple $(\psi,\pi',\alpha)$ is valid, and this shows that $\pi\le\psi$.
\endproof

\proof{Proof of Theorem~\ref{thm:non-dominated}} ($\Leftarrow$) Theorem 6.34 in~\cite{conforti2014integer} shows that if conditions $(a)$--$(e)$ are satisfied, then $(\psi,\pi,\alpha)$ is a minimal valid tuple for $M_b$. 
Since $\alpha = 1$, the tuple is nontrivial.\medskip

($\Rightarrow$) Suppose that $(\psi, \pi, \alpha)$ is a nontrivial minimal valid tuple for $M_b$.

\begin{enumerate}[$(a)$]
\item  This proof is the same as the subadditivity proof in Proposition~\ref{prop:psi_sublinear}.

\item  We first establish the following claim.

\begin{claim}
$\psi(r) \geq \sup_{\epsilon > 0}\frac{\pi(\epsilon r)}{\epsilon} = \lim_{\epsilon \to 0^+}\frac{\pi(\epsilon r)}{\epsilon} = \lim\sup_{\epsilon \to 0^+}\frac{\pi(\epsilon r)}{\epsilon}$.
\end{claim}

\begin{claimproof}
Since $(\psi, \pi, \alpha)$ is minimal, from Proposition~\ref{prop:psi_sublinear} $\psi$ is sublinear and $\pi\leq \psi$. Hence for $\epsilon>0$ and $r\in \R^n$ we have that $\frac{\pi(\epsilon r)}{\epsilon}\leq \frac{\psi(\epsilon r)}{\epsilon} = \psi(r)$.
Thus, $\sup_{\epsilon > 0}\frac{\pi(\epsilon r)}{\epsilon} \leq \psi(r)$ and this implies that $\sup_{\epsilon > 0}\frac{\pi(\epsilon r)}{\epsilon}$ is a finite real number. By Theorem 7.11.1 in \cite{hillePhillips} and the subadditivity of $\pi$, this implies that $\sup_{\epsilon > 0}\frac{\pi(\epsilon r)}{\epsilon} = \lim_{\epsilon \to 0^+}\frac{\pi(\epsilon r)}{\epsilon} = \lim\sup_{\epsilon \to 0^+}\frac{\pi(\epsilon r)}{\epsilon}$.\end{claimproof}

The above claim shows that the function $\psi'(r) := \lim_{\epsilon \to 0^+} \frac{\pi(\epsilon r)}{\epsilon}$ is well defined, and $\psi' \leq \psi$. Furthermore, by Lemma \ref{lem:psi-using-sup}, $\psi'$ is sublinear. We prove below that $(\psi', \pi, \alpha)$ is a valid tuple. Therefore, since $(\psi, \pi, \alpha)$ is minimal, validity of $(\psi', \pi, \alpha)$ will imply that $\psi = \psi'.$

Assume by contradiction that $(\psi', \pi, \alpha)$ is not valid. Then there exists $(s,y)\in M_b$ such that
$$\sum_{r\in \R^n}\psi'(r)s(r)+\sum_{p\in \R^n}\pi(p)y(p) = \alpha - \delta$$
for some $\delta >0$. Define $\tilde{r} := \sum_{r\in \R^n} r s(r)$. Since $\psi'(r) = \lim_{\epsilon \to 0^+} \frac{\pi(\epsilon r)}{\epsilon}$, there exists some $ a>0$ such that
\begin{equation}\label{eq:delta}
\frac{\pi(\epsilon \tilde{r})}{\epsilon} < \psi'(\tilde{r})+\delta~~\text{for all}~~0<\epsilon<  a.
\end{equation}

Let $D\in \Z_{>0}$ be such that $1/D<a$ and define $\tilde{y}$ to be
$$\tilde{y}(r) :=
\begin{cases}
y(r)+D&\text{if }r = \tilde{r}/D,\\
y(r)&\text{if }r \neq \tilde{r}/D.
\end{cases}$$
Note that
$$\sum_{r\in \R^n}r\tilde{y}(r) = \sum_{r\in \R^n}rs(r)+\sum_{p\in \R^n}py(p) \in b+\Z^n,$$
and so $(0,\tilde{y})\in M_b$. Hence $\sum_{p\in \R^n}\pi(p)\tilde{y}_p \geq \alpha$. However,
\begin{align*}
\sum_{p\in \R^n}\pi(p)\tilde{y}(p) & = \frac{\pi(\tilde{r}/D)}{1/D}+\sum_{p\in \R^n}\pi(p)y(p) \\
&< \psi'(\tilde{r})+\delta + \sum_{p\in \R^n}\pi(p)y(p)&&\text{by \eqref{eq:delta} and by definition of $D$}\\
& \leq  \sum_{r\in \R^n}\psi'(r)s(r) +\delta + \sum_{p\in \R^n}\pi(p)y(p) &&\text{by sublinearity of $\psi'$}\\
& = \alpha,
\end{align*}
which is a contradiction.

\item We now show that $\pi$ is Lipschitz continuous with Lipschitz constant $L:=\max_{\|r\|=1}\psi(r)$. By Proposition~\ref{prop:psi_sublinear}, $\psi$ is sublinear; thus, it is continuous. Therefore $\max_{\|r\|=1}\psi(r)$ is attained. Moreover, by subadditivity of $\pi$, we obtain that $\pi(x) - \pi(y) \leq \pi(x - y)$ for all $x,y \in \R^n$. Therefore, $|\pi(x) - \pi(y)| \leq \max\{\pi(x - y),\pi(y-x)\}$. Thus, for all $x\neq y$,
$$\frac{|\pi(x)-\pi(y)|}{\|x-y\|} \leq \frac{\max\{\pi(x-y),\pi(y-x)\}}{\|x-y\|} \leq \frac{\max\{\psi(x-y),\psi(y-x)\}}{\|x-y\|} \leq L,$$ where the second inequality follows from Proposition~\ref{prop:psi_sublinear}, and the last inequality from the positive homogeneity of $\psi$.

\item We prove this with a sequence of claims.

\begin{claim}\label{claim:nonnegative_rational}
$\pi(r) \geq 0$ for all $r\in \R^n$.
\end{claim}

\begin{claimproof}
Let $p^*\in \Q^n$. Then there exists $D\in \Z_{>0}$ such that $Dp^*\in \Z^n$. Let $(s,y)\in M_b$ and, for some $k\in \Z_+$, define $(s,\tilde{y})$ by setting $\tilde{y}({p^*}) := y({p^*})+kD$ and $\tilde{y}(p) := y(p)$ for $p\ne p^*$. Note that $(s, \tilde{y})\in M_b$ for every $k\in\Z_+$. This shows that $\pi(p^*)\ge0$ for every $p^*\in\Q^n$.
Since $\pi$ is Lipschitz continuous by part $(c)$ above, we must have $\pi \geq 0$ everywhere.
\end{claimproof}

\begin{claim}
$\pi(z) = 0$ for all $z\in \Z^n$.
\end{claim}

\begin{claimproof}
Assume to the contrary that there is some $z\in \Z^n$ such that $\pi(z) \neq 0$. By the previous claim, $\pi(z)>0$. Define $\pi'$ to be $\pi'(z):=0$ and $\pi'(p):=\pi(p)$ for $p\neq z$. Then $(\psi,\pi',\alpha)$ is easily seen to be a valid tuple. This contradicts the minimality of $(\psi, \pi, \alpha)$. 
\end{claimproof}

We now show that $\alpha = 1$. Since $\psi, \pi \geq 0$ by part $(b)$ and Claim \ref{claim:nonnegative_rational}, if $\alpha = 0$ or $\alpha=-1$ then this would contradict the fact that the tuple is nontrivial.

\item The proof is identical to part $(d)$ of the proof of Theorem 6.22 in~\cite{conforti2014integer}.

\end{enumerate}

This concludes the proof of the theorem.
\endproof

\section{Continuity of linear functional $\hat F_\pi$.}

Just like $F_{\psi, \pi}$, define the linear functional $\hat F_\pi$ on $\R^{(\R^n)}$ by $$\hat F_\pi (y):= \sum_{p\in \R^n} \pi(p)y(p).$$

\begin{lemma}\label{lem:cts-dir-der} Let $\pi :\R^n \to \R$.  If $\hat F_\pi$ is a continuous linear functional with respect to the $|\cdot|_*$ norm, then $\sup_{\epsilon > 0}\frac{\pi(\epsilon r)}{\epsilon} < \infty$ for all $r \in \R^n\setminus\{0\}$.


\end{lemma}

\proof{Proof.} $\hat F_\pi$ is a continuous linear functional if and only if it is bounded, i.e., there exists $M\in \R$ such that $\frac{|\hat F_\pi(y)|}{|y|_*} \leq M$ for all $y \neq 0$ (see Conway~\cite[Chapter III, Proposition 2.1]{conway2013course}). 

Suppose $\hat F_\pi$ is bounded and suppose to the contrary that for some $\bar r \in \R^n\setminus\{0\}$, $\sup_{\epsilon > 0}\frac{\pi(\epsilon \bar r)}{\epsilon} = \infty$. For every $\epsilon > 0$, consider $y_\epsilon \in \R^{(\R^n)}$ given by $$y_\epsilon(r) =  \left\{\begin{array}{cl}1 & \textrm{if } r = \epsilon\bar r, \\ 0 & \textrm{otherwise.} \end{array}\right.$$

Observe that $$\sup_{\epsilon > 0}\frac{|\hat F_\pi(y_\epsilon)|}{|y_\epsilon|_*} = \sup_{\epsilon>0}\frac{|\pi(\epsilon \bar r)|}{\|\epsilon\bar r\|} = \frac{1}{\|\bar r\|}\sup_{\epsilon>0}\frac{|\pi(\epsilon \bar r)|}{\epsilon} = \infty.$$
This shows that $\hat F_\pi$ is not bounded, contradicting the assumption.
\endproof

\bibliographystyle{plain}
\bibliography{../full-bib}

\end{document}